\documentclass[12pt,reqno]{amsart}
  \usepackage[all,2cell,arrow]{xy}
  \usepackage{amsfonts}
  \usepackage{amsthm}
  \usepackage{amsmath}
  \usepackage{amssymb}
\usepackage{mathtools}
\usepackage{tensor}
 \numberwithin{equation}{section}
 \usepackage{color}
\usepackage{graphicx}
  \usepackage{hyperref}
\usepackage{enumerate,xspace}

\UseAllTwocells
\CompileMatrices

\renewcommand{\epsilon}{\varepsilon}
\renewcommand{\phi}{\varphi}

\newcommand{\ca}{\ensuremath{\mathcal A}\xspace}
\newcommand{\cb}{\ensuremath{\mathcal B}\xspace}
\newcommand{\cc}{\ensuremath{\mathcal C}\xspace}
\newcommand{\cd}{\ensuremath{\mathcal D}\xspace}
\newcommand{\ce}{\ensuremath{\mathcal E}\xspace}
\newcommand{\cf}{\ensuremath{\mathcal F}\xspace}

\newcommand{\cp}{\ensuremath{\mathcal P}\xspace}

\newcommand{\cs}{\ensuremath{\mathcal S}\xspace}

\newcommand{\cv}{\ensuremath{\mathcal V}\xspace}

\newcommand{\cx}{\ensuremath{\mathcal X}\xspace}

\newcommand{\bbc}{\ensuremath{\mathbb C}\xspace}

\newcommand{\bbn}{\ensuremath{\mathbb N}\xspace}
\newcommand{\bbp}{\ensuremath{\mathbb P}\xspace}

\newcommand{\Span}{\ensuremath{\mathbf{Span}}\xspace}

\newcommand{\Set}{\ensuremath{\mathbf{Set}}\xspace}
\newcommand{\OpCat}{\ensuremath{\mathbf{OpCat}}\xspace}
\newcommand{\SkewMon}{\ensuremath{\mathbf{SkewMonCat}}\xspace}

\newcommand{\Ord}{\ensuremath{\mathbf{Ord}}\xspace}

\newcommand{\op}{\ensuremath{^{\textnormal{op}}}}

\DeclareMathOperator{\cod}{cod}
\DeclareMathOperator{\dom}{dom}

\DeclareMathOperator{\Coll}{Coll}
\DeclareMathOperator{\sColl}{Coll}

\DeclareMathOperator{\Bq}{Bq}
\DeclareMathOperator{\Op}{Op}
\DeclareMathOperator{\Epi}{Epi}

\newcommand{\ox}{\otimes}
\newcommand{\x}{\times}

\def\1c#1{\stackrel{#1}{\to}}

  \newtheorem{proposition}{Proposition}[section]
  \newtheorem{lemma}[proposition]{Lemma}

  \newtheorem{theorem}[proposition]{Theorem}

  \theoremstyle{definition}
  \newtheorem{definition}[proposition]{Definition}
  \newtheorem{notation}[proposition]{Notation}
  \newtheorem{example}[proposition]{Example}

  \theoremstyle{remark}
  \newtheorem{remark}[proposition]{Remark}

  \newcounter{c}
  \renewcommand{\[}{\setcounter{c}{1}$$}
  \newcommand{\etyk}[1]{\vspace{-7.4mm}$$\begin{equation}\Label{#1}
  \addtocounter{c}{1}}
  \renewcommand{\]}{\ifnum \value{c}=1 $$\else \end{equation}\fi}
\begin{document}

 \title[Operadic categories]{Operadic categories and their skew monoidal categories of collections}

\author{Stephen Lack}
\address{Department of Mathematics, Macquarie University NSW 2109, Australia}
\email{steve.lack@mq.edu.au}

\date{\today}

\begin{abstract}
I describe a generalization of the notion of operadic category due to Batanin and Markl. For each such operadic category I describe a skew monoidal category of collections, such that a monoid in this skew monoidal category is precisely an operad over the operadic category. In fact I describe {\em two} skew monoidal categories with this property. The first has the feature that the operadic category can be recovered from the skew monoidal category of collections; the second has the feature that the right unit constraint is invertible. In the case of the operadic category \cs of finite sets and functions, for which an operad is just a symmetric operad in the usual sense, the first skew monoidal category has underlying category $[\bbn,\Set]$, and the second is the usual monoidal category of collections $[\bbp,\Set]$ with the substitution monoidal structure. 
\end{abstract}
  
\maketitle


\section{Introduction}

In the beginning \cite{May-operads}, an operad was a formalism for describing certain
sorts of internal structures in a symmetric monoidal category. For
each natural number $n$ one specified the $n$-ary operations which
could be defined in the structure as well as various equations which
hold between these operations. There would typically be an action of the
symmetric group $S_n$ which allowed the input variables to be
permuted, but there was also a ``non-symmetric'' or ``plain'' flavour
of operad which did not involve these actions.

The renaissance of operads \cite{renaissance} which was celebrated in the mid-1990s saw
not just a renewed interest in operads but an explosion of new flavours
of operad. These included ``coloured'' versions (also known as
multicategories, symmetric or otherwise), higher globular operads \cite{Batanin-MonoidalGlobular}, and modular operads \cite{GetzlerKapranov-Modular}. There were
various approaches to incorporate many of these into a single
framework \cite{Hermida-CoherentUniversal,Leinster-book}

This expansion in the scope of operads has continued, and several
new frameworks have appeared recently. One of these is the {\em operator
  categories} of Barwick \cite{Barwick-Operator}; another, more
general, is the {\em operadic categories} of the title, introduced by
Batanin and Markl in \cite{BataninMarkl-Operadic}.  For each such operadic category, there is a corresponding notion of operad. Thus there is an operadic category for symmetric operads, another for plain operads, another for coloured versions of these (for a given set of colours), and still others for other notions of operad. These operadic categories were put to spectacular use in proving a duoidal version of the Deligne conjecture. 

In the original (symmetric) operads, the category \bbp of finite cardinals and bijections (equivalently, the disjoint union of the symmetric groups) plays a key role. For a symmetric monoidal category \cv with colimits preserved by tensoring, the functor category $[\bbp,\cv]$ has a (non-symmetric) monoidal structure, and a monoid with respect to this monoidal structure is precisely a symmetric operad in \cv \cite{Kelly-operads}. An object of $[\bbp,\cv]$ is sometimes called a {\em collection} in \cv, and consists of an object $T_n$ of \cv for each $n\in\bbn$, equipped with an action of the symmetric group. The monoidal category of such collections is sometimes written as $\Coll(\cv)$. 

In the case of plain operads there are no actions of the symmetric groups, and so rather than \bbp one uses the discrete category \bbn consisting of just the finite cardinals and identity morphisms; there is once again a monoidal structure on $[\bbn,\cv]$ with respect to which the monoids are the plain operads in \cv. 

In this paper, I shall introduce a mild generalization of the operadic
categories of \cite{BataninMarkl-Operadic}, and for each such
``generalized operadic category'' \cc I shall define a {\em skew}
monoidal category \cite{Szlachanyi-skew, skew} of collections
$\sColl_\cc(\cv)$. Skew monoidal categories are similar to monoidal
categories except that the unit and associativity maps are not
required to be invertible. The most important case of this
construction is where \cv is just the category \Set of sets and functions, equipped with the usual cartesian monoidal structure. The skew monoidal category $\sColl_\cc(\Set)$ is equipped with an opmonoidal functor into  $[\bbn,\Set]$. I shall write $\sColl_\cc$ for $\sColl_\cc(\Set)$. 

I shall show that the generalized operadic category \cc can be recovered from
$\sColl_\cc$ along with its opmonoidal functor into $[\bbn,\Set]$, and
I characterize which skew monoidal categories over $[\bbn,\Set]$ arise
in this way, and further characterize those corresponding to the
genuine operadic categories of \cite{BataninMarkl-Operadic}. This
provides a new, but equivalent, definition of operadic category, as
well as a rather different point of view. I
regard this as the main contribution of the paper. (The mild
generalization in the 
definition seems far less important, although it does allow a cleaner
way for presheaves to be seen as operads.)

If we start with the operadic category \cs for symmetric operads, the
resulting skew monoidal category $\sColl_\cs$ is {\em not} just
$\Coll(\Set)=[\bbp,\Set]$; in fact as a category it is $[\bbn,\Set]$,
but with a different skew monoidal structure to that mentioned
above. Nonetheless, there is a way to recover $\Coll(\Set)$ from
$\sColl_\cs$. For a large class of skew monoidal categories
$(\ce,*,U)$, there is a way \cite{mw} to associate a new skew monoidal
category $\ce^U$ for which the right unit map $X\to X\ox U$ is
invertible, and such that the two skew monoidal categories have the
same category of monoids. (A skew monoidal category for which the
right unit maps are invertible is said to be {\em right normal}.) When this construction is applied to $\sColl_\cs$ the resulting skew
monoid\-al category $\sColl^U_\cs$ is in fact monoidal, and is
monoidally equivalent to $\Coll(\Set)$. This is the second main
contribution of the paper.

When the construction is applied to $\sColl_\cc$ for a general operadic category \cc, the resulting skew monoidal category $\sColl^U_\cc$ may not be a monoidal category, but I give sufficient conditions under which it is so. 

The structure of the paper is as follows. I begin, in
Section~\ref{sect:operadic}, with the generalized notion of operadic
category, its relationship to the operadic categories of
\cite{BataninMarkl-Operadic}, and a few key examples. In fact I use
``operadic category'' for the new more general notion, and speak of
``genuine operadic categories'' when I wish to refer to the original
notion of \cite{BataninMarkl-Operadic}.  In
Section~\ref{sect:sColl}, I define the skew monoidal category
$\sColl_\cc$ associated to an operadic
category. Section~\ref{sect:operads} is about \cc-operads for an operadic
category \cc; that is, about monoids in $\sColl_\cc$. I describe in Section~\ref{sect:functoriality}  the
dependence of the $\sColl_\cc$ construction on \cc, and in
Section~\ref{sect:V} the modifications needed when collections are 
taken in a monoidal category \cv other than $\Set$. Then in 
Section~\ref{sect:characterization} the characterization of those skew monoidal
categories over $\sColl_\cs$ arising from an operadic category is given. In
Sections~\ref{sect:bbc} and~\ref{sect:normalization} I show how to
replace $\sColl_\cc$ with a right normal skew monoidal category
$\sColl^U_\cc$ with the same category of monoids (that is, operads),
as well as giving a sufficient condition for $\sColl^U_\cc$ to be
monoidal. Finally in Section~\ref{sect:examples} I
describe various examples of operadic categories, and determine in
each case whether the sufficient condition holds.

\subsection*{Acknowledgements.} Most of the material in this paper was
presented in an invited talk at the CT2014 conference in Cambridge. I
am grateful to the organizers for their generous invitation, for the
opportunity to speak, and for a very enjoyable conference. It is also
a pleasure to acknowledge various helpful comments from Michael
Batanin and other members of the Australian Category Seminar. The MSc
thesis \cite{Andrianopoulos-thesis} of Jim Andrianopolous was an important
influence. Finally I gratefully acknowledge the support of an Australian Research Council Discovery Project DP130101969 and a Future Fellowship FT110100385.

\section{Operadic categories} \label{sect:operadic}

For an object $d$ of a category \cc, the ``slice category'' $\cc/d$ has morphisms with codomain $d$ as objects, and commutative triangles as morphisms. It is equipped with a functor $\dom\colon \cc/d\to \cc$ sending an object of $\cc/d$ to the domain of the corresponding morphism.
\[ \xymatrix @R1pc @C1pc {
& \cc/d \ar[rrrr]^{\dom} &&&& \cc \\
b \ar[rr]^{\phi} \ar[dr]_{\psi\phi} && c \ar[dl]^{\psi} & \mapsto & b \ar[rr]^{\phi} && c \\
& d } \] 
For any functor $F\colon \cx\to \cc$ and any object $x\in\cx$ there is an induced functor $F/x\colon \cx/x\to \cc/Fx$ sending $\psi\colon y\to x$ to $F\psi\colon Fy\to Fx$. 
In particular, for the functor $\dom\colon \cc/d\to \cc$ and any $\psi\colon c\to d$, the induced $\dom\!/\psi\colon (\cc/d)/(\psi\colon c\to d)\to \cc/c$ is an isomorphism of categories. 

For a set $I$, there are two possible ways to define the category of $I$-indexed sets: as $\Set^I$ or as $\Set/I$. Of course these two are equivalent, via the functor $\Set/I\to\Set^I$ sending a set over $I$ to its fibres. 

\begin{notation}
  Throughout this paper \cs will denote (any skeleton of) the category of finite sets.
\end{notation}

If $I$ is a finite set (in \cs) then the equivalence $\Set/I\simeq \Set^I$ clearly restricts to an equivalence $R_I\colon \cs/I\to \cs^I$. This is determined only up to isomorphism; we shall suppose a fixed choice to have been made.

If we move from categories of sets to some other category \cc, these
two approaches to families are no longer equivalent, or even directly
comparable. There is still the category $\cc^I$ of $I$-indexed
families in \cc: this is an ``external'' notion of family. But there
is also the ``internal'' version of indexed family, where the
indexation is done using an object $X\in\cc$ rather than a set $I$,
and now $\cc/X$ can be thought of as the category of ``$X$-indexed
families of objects in \cc''. Both of these are important: among other
things, the first is fundamental to the theory of enriched categories and the second is fundamental to the theory of internal categories. 
In an operadic category the equivalence between the internal and
external notions is partially restored. 

I shall now introduce the
promised generalization of the operadic categories of
\cite{BataninMarkl-Operadic}. For the precise relationship between the two
definitions, see Proposition~\ref{prop:comparison} below. 

\begin{definition}\label{defn:operadic}
  An {\em operadic category} is a category \cc equipped with the following structure:
  \begin{enumerate}
  \item a functor $|~|\colon \cc\to\cs$, which we call the {\em cardinality functor};
\item for each object $c\in\cc$ a functor $R_c\colon \cc/c\to \cc^{|c|}$ making the diagram 
\[\xymatrix{ \cc/c \ar[r]^{R_c} \ar[d]_{|~|} & \cc^{|c|} \ar[d]^{|~|} \\
\cs/|c| \ar[r]_{R_{|c|}} & \cs^{|c|} } \]
commute;
  \end{enumerate}
subject to some conditions to which we shall soon turn. First, however, we introduce some notation and terminology. For a morphism $\psi\colon c\to d$ in \cc, the functor $R_d$ gives a $|d|$-indexed family of objects of \cc. The $i$th of these, for some $i\in |d|$, will be written as $\psi^{-1}i$; these $\psi^{-1}i$ will be called the {\em fibres} of $\psi$. For morphisms $\phi\colon b\to c$ and $\psi\colon c\to d$, seen as defining a morphism $\phi\colon (\psi\phi)\to \psi$ in $\cc/d$, we sometimes write $\phi^\psi\colon R(\psi\phi)\to R(\psi)$ for its image under $R_d$.

An object $u\in\cc$ is said to be {\em trivial} if $|u|=1$ and $R_u=\dom$.  
The commutativity of the square implies that $|\phi^{-1}i|=|\phi|^{-1}i$. In particular, if $u=1^{-1}_c i$ is a fibre of an identity morphism, then $|u|=1$.
We shall often omit the subscript and simply write $R$ for $R_d$. 

We now turn to the conditions.
\begin{enumerate}\addtocounter{enumi}{2}
\item Any fibre $1^{-1}_ci$ of an identity morphism is trivial;
\item For any morphism $\psi\colon c\to d$, the diagram
\[ \xymatrix{
(\cc/d)/(\psi\colon c\to d) \ar[r]^-{R/\psi} \ar[d]_{\dom\!/\psi} & 
\cc^{|d|}/R\psi \ar@{=}[r]   & 
{}\prod\limits_{j\in|d|} \cc/\psi^{-1}j \ar[d]^{\Pi_j R} \\
\cc/c \ar[r]_{R} & \cc^{|c|} \ar@{=}[r] & {}\prod\limits_{j\in|d|} \cc^{|\psi^{-1}j|} } \]
commutes (the ``double slice condition''). Here the ``equality'' on
the lower line is defined using the equality
$|\psi^{-1}j|=|\psi|^{-1}j$ and the canonical isomorphism
$|c|\cong\sum_j |\psi|^{-1}j$. The object part of this double slice
condition says that for a composable pair $\phi\colon c\to d,
\psi\colon d\to e$, and $i\in|d|$, the equation
$(\phi^\psi_{|\psi|j})^{-1}i=\phi^{-1}i$ holds. 
\end{enumerate}
We often write $U$ for the set of all trivial objects. 
\end{definition}

\begin{proposition}\label{prop:trivial}
  In an operadic category, an object $u$ is trivial if and only if it is a fibre of $1_u$.
\end{proposition}

\proof
Any fibre of an identity morphism is trivial. Conversely, if $u$ is trivial, then $|u|=1$ and the (unique) fibre of any morphism $c\to u$ is $c$; in particular, the fibre of $1_u$ is $u$.
\endproof

We record the relationship with the operadic categories of Batanin-Markl as the following proposition. In fact the only notion defined in \cite{BataninMarkl-Operadic} is called a {\em strict} operadic category; we shall call it a {\em genuine} operadic category when we wish to distinguish it from the ``generalized'' operadic categories considered here. (We never deal with the ``non-strict'' notion, in which the commutative diagrams in conditions (2) and (4) of the definition are replaced by isomorphisms, satisfying as yet unspecified coherence conditions.)

\begin{proposition}\label{prop:comparison}
  A strict operadic category in the sense of \cite{BataninMarkl-Operadic} is precisely an operadic category, in the sense of Definition~\ref{defn:operadic}, in which each connected component has a chosen terminal object, these objects are trivial, and they are the only trivial objects. 
\end{proposition}

A {\em strict operadic functor} between operadic categories \cc and \cd will be a functor $F\colon\cc\to\cd$ which strictly commutes with both the functors into \cs and the functors $R$, in the sense that the diagrams
\[\xymatrix{
\cc \ar[rr]^{F} \ar[dr]_{|~|} && \cd \ar[dl]^{|~|} && \cc/c \ar[r]^{F/c} \ar[d]_{R_c} & \cd/Fc \ar[d]^{R_{Fc}} \\
& \cs & && \cc^{|c|} \ar[r]_{F^{|c|}} & \cd^{|Fc|} }\]
commute; the second makes sense because $|Fc|=|c|$ by commutativity of the first. 
We write \OpCat for the category of operadic categories and strict operadic functors. 

\begin{proposition}
A strict operadic functor sends trivial objects to trivial objects. 
\end{proposition}

\proof
If $F\colon \cc\to\cd$ is a strict operadic functor and $u$ is trivial in \cc, then $|Fu|=|u|=1$; while $u$ is a fibre of $1_u$ and so $Fu$ is a fibre of $F1_u=1_{Fu}$; thus $Fu$ is also trivial. 
\endproof

\begin{remark}
  The strict operadic functors of \cite{BataninMarkl-Operadic} are required to strictly preserve the chosen terminal objects, but this just amounts to preserving the trivial objects. 
\end{remark}

\begin{example}
  The category \cs itself is operadic, with $|~|$ given by the
  identity functor. In fact \cs is the terminal operadic category, in
  the sense that for any operadic category \cc, there is a unique strict operadic functor $\cc\to\cs$, given by the cardinality functor $|~|$. 
\end{example}

\begin{example}
  The category \cp of finite ordinals and order-preserving functors is
  operadic; the cardinality functor forgets the order, and the $R$ are
  constructed using the fibres with their induced ordering. (The
  category \cp has been given various names over the years. It contains
  the simplex category $\mathbf{\Delta}$ as the full subcategory of
  all non-empty finite ordinals. When \cp is made into a monoidal
  category via ordinal sum, it is sometimes called the ``algebraists'
  $\Delta$'', but this monoidal structure will not be used here. The
  letter $\cp$ has been chosen to suggest {\em plain} operads, as
  opposed to the {\em symmetric} operads corresponding to \cs. 
\end{example}

\begin{example}
As observed in \cite{BataninMarkl-Operadic}, a {\em category of operators} in the sense of Barwick \cite{Barwick-Operator} is an operadic category \cc with finite hom-sets and a  terminal object $1$, in which the cardinality functor is the \cs-valued representable functor $\cc(1,-)$, and the fibres are the actual fibres defined using pullback.  
\end{example}

Our first example which does {\em not} satisfy the extra condition in \cite{BataninMarkl-Operadic} is:

\begin{example}\label{ex:x}
Any category \ca can be made into an operadic category by defining the cardinality $|a|$ of any object $a$ to be $0$. Then there are no fibres and no trivial objects. 
\end{example}

Batanin and Markl also describe how to make any category \ca into an
operadic category in their sense: one freely adds a terminal object
and makes this object trivial, while all objects from the original
category have cardinality $0$. I shall write $\ca_1$ for this operadic
category, and I shall have more to say about the difference between
$\ca$ and $\ca_1$ below. The next example is closely related to \cite[Chapter~3]{Andrianopoulos-thesis}. 

\begin{example}\label{ex:jim}
Any category \ca can be made into an operadic category by defining the
cardinality $|a|$ of any object to be $1$, and defining each $R_a$ to
be the domain functor.   
\end{example}

There are many further examples of operadic category given in \cite{BataninMarkl-Operadic}; some of these are discussed in Section~\ref{sect:examples} below.

\section{The skew monoidal category of collections} \label{sect:sColl}

A skew monoidal category \cite{Szlachanyi-skew, skew} is a category \ce equipped with a functor $\ce\x\ce\to\ce$, whose effect on $(X,Y)$ is written $X*Y$, an object $U$, and natural transformations
\[ \xymatrix @R0.5pc {
(X*Y)*Z \ar[r]^{\alpha} & X*(Y*Z) \\
U*X \ar[r]^{\lambda} & X \\
X \ar[r]^{\rho} & X*U } \]
subject to five axioms which are recalled below. These natural transformations are not required to be invertible, but it is useful to be able to discuss the case when some or all of them are so. The skew monoidal category is said to be {\em left normal} if $\lambda$ is invertible and {\em right normal} if $\rho$ is invertible. It is said to be {\em Hopf} if $\alpha$ is invertible; of course if all three are invertible then it is just a monoidal category. 

In this section we shall show how to construct a skew monoidal category $\sColl_\cc(\Set)$, or $\sColl_\cc$ for short, from any operadic category \cc. 

Let \cc be an operadic category, and write $C$ for the set of objects of \cc. 
The underlying category of $\sColl_\cc$ will be the slice category $\Set/C$.
An object of $\Set/C$ consists of a set $X$ equipped with a function $\partial\colon X\to C$, but we normally regard $\partial$ as understood and simply call the object $X$.
Given such an object $X$ and an element $c\in C$ we write $X_c$ for
the (actual!) fibre $\partial^{-1}(c)$. For $x\in X$ we write $|x|$
for the set $|\partial(x)|$.

The tensor product $X*Y$ of $X$ and $Y$ is given by the formula
$$(X*Y)_c = \sum_{\phi\colon c\to d} X_d\x \prod_{i\in|d|} Y_{\phi^{-1}i}.$$
Thus an element of $X*Y$ consists of a morphism $\phi\colon c\to d$ in \cc, an element $x\in X_d$, and a $|d|$-indexed family $y=(y_i)_{i\in|d|}$ with $y_i\in Y_{\phi^{-1}i}$.
The function $\partial\colon X*Y\to C$ sends such an element $(x,\phi,y)$ to the domain of $\phi$.
This clearly extends to a functor $\Set/C\x \Set/C\to \Set/C$ sending $(X,Y)$ to $X*Y$.

The unit is the set $U$ of trivial objects of \cc, with $\partial$ given by the inclusion $U\to C$.

The remaining structure in a skew monoidal category consists of the natural transformations $\lambda$, $\rho$, and $\alpha$, to which we now turn.
An element of $(U*X)_c$ has the form $(u,\phi,x)$, where $\phi\colon c\to u$ and $u\in U$; since $|u|=1$ and the unique fibre of $\phi$ is $c$, $x$ just consists of a single element of $X_c$. We may now define {\em the left unit constraint} for our skew monoidal structure to be 
\[ \xymatrix @R0pc {
U*X \ar[r]^{\lambda} & X \\
(u,\phi,x) \ar@{|->}[r] & x } \]  
which is clearly natural in $X$.

\begin{remark}\label{rmk:lambda}
  The left unit map $\lambda$  is invertible if and only if every object $c$ has a unique map to some trivial object $u$. If this is the case, we may write $!_c\colon c\to u(c)$ for this map. For any morphism $\phi\colon c\to d$, the composite 
\[ \xymatrix{ c \ar[r]^{\phi} & d \ar[r]^{!_d} & u(d) } \] 
is a morphism to a trivial object, so by uniqueness $u(d)=u(c)$ and this composite is the unique map. Thus it follows that each connected component of \cc has a chosen terminal object, and this terminal object is trivial. It further follows that these chosen terminal objects are the only trivial objects. Thus $\lambda$ will be invertible if and only if our ``generalized'' operadic category is a genuine operadic category in the sense of \cite{BataninMarkl-Operadic}.
\end{remark}

We shall say that a morphism $\phi$ is {\em fibrewise trivial} if all of its fibres are trivial. Part of the definition of operadic category is that identity morphisms are fibrewise trivial. 
In the operadic categories \cp and \cs a morphism is fibrewise trivial if and only if it is bijective. 

An element of $(X*U)_c$ has the form $(x,\phi,u)$, where $\phi\colon
c\to \partial(x)$ must be fibrewise trivial in order to define $u$. Thus we may define {\em the right unit constraint} for our skew monoidal structure to be
\[ \xymatrix @R0pc {
X \ar[r]^{\rho} & X*U \\
x \ar@{|->}[r] & (x,1_{\partial(x)},R1_{\partial(x)}) } \]
which once again is clearly natural in $X$.

\begin{remark}\label{rmk:rho}
  The right unit map $\rho$ will be invertible if and only if the only
  fibrewise trivial morphisms are the identities. This is true in \cp
  but not in \cs (or most other examples).
\end{remark}

Next we turn to the associativity map $\alpha$. First we describe $(X*Y)*Z$ and $X*(Y*Z)$. An element of $((X*Y)*Z)_c$ involves a morphism $\phi\colon c\to d$, an element of $(X*Y)_d$, and a $|d|$-indexed family $z$ with $z_i\in Z_{\phi^{-1}i}$; and an element of $(X*Y)_d$ will consist of a morphism $\psi\colon d\to e$, an element $x\in X_e$, and a $|e|$-indexed family $y$ with $y_i\in Y_{\psi^{-1}i}$. We denote such an object with $(x,\psi,y,\phi,z)$. 

An element of $(X*(Y*Z))_c$ consists of a morphism $\theta\colon c\to e$, an element $x\in X_e$, and an $|e|$-indexed family $(y,\tau,z)$ with $(y,\tau,z)_j\in (Y*Z)_{\theta^{-1}j}$.
Here $\tau_j\colon \theta^{-1}j\to v_j$ and $y_j\in Y_{v_j}$, while $z_j$ is a $|v_j|$-indexed family with $(z_j)_i\in Z_{\tau^{-1}_ji}$. We may collect all the $v_j$ into an object $v\in \cc^{|e|}$, and all the $\tau_j$ into a single morphism $\tau\colon R(\theta)\to v$ in $\cc^{|e|}$. We now define {\em the associativity constraint} for the skew monoidal structure to be 
\[ \xymatrix @R0pc {
(X*Y)*Z \ar[r]^{\alpha} & X*(Y*Z) \\
(x,\psi,y,\phi,z) \ar@{|->}[r] & (x,\psi\phi,y,\phi^\psi,z) } \]
which is once again clearly natural. (Recall that $\phi^\psi\colon
R(\psi\phi)\to R(\psi)$ is the image under $R_d\colon \cc/d\to
\cc^{|d|}$ of the morphism $\phi\colon (\psi\phi\colon b\to d)\to
(\psi\colon c\to d)$ in $\cc/d$.)

\begin{remark}\label{rmk:alpha}
When is $\alpha$ invertible? In particular this would imply that 
$\alpha\colon (C*C)*C\to C*(C*C)$ is invertible. Given $\theta\colon c\to e$ and
$\tau\colon R(\theta)\to d$, in the form $\tau_j\colon \theta^{-1}j\to
d_j$ for $j\in |e|$, there needs to be a unique $\phi\colon c\to d$
and $\psi\colon d\to e$ with $\psi\phi=\theta$ and
$\phi^\psi_j=\tau_j$; that is, a unique $\phi\colon \theta\to\psi$ in
$\cc/e$ with $R\psi=d$ and $R(\phi)=\tau$.  This in turn says that each $R\colon \cc/e\to \cc^{|e|}$ is a discrete opfibration. Conversely, it is not hard to check that in this case not just $\alpha\colon (C*C)*C\to C*(C*C)$, but all the components of $\alpha$ are invertible.

In the example of \cs, the $R$ functors are equivalences, but are not discrete opfibrations, so $\alpha$ is not invertible. In the case of \cp, however, the $R$ are in fact isomorphisms, so $\alpha$ is invertible.
\end{remark}

Now we turn to the axioms \cite{Szlachanyi-skew, skew} for skew monoidal categories. For monoidal categories (where $\alpha$, $\lambda$, and $\rho$ are all invertible) two axioms suffice, but for skew monoidal categories five are needed.

\subsection*{The $(\lambda,\rho)$-compatibility condition.}
This says that the composite 
$$\xymatrix{
U \ar[r]^-{\rho} & U*U \ar[r]^-{\lambda} & U }$$
is the identity. Since each $U_c$ has at most one element, this is obviously true. 
 
\subsection*{The $(\alpha,\lambda)$-compatibility condition.}
This says that the diagram
$$\xymatrix{
(U*X)*Y \ar[r]^{\alpha} \ar[dr]_{\lambda*1} & U*(X*Y) \ar[d]^{\lambda}
\\ & X*Y }$$
commutes. Now
\begin{align*}
\lambda(\alpha(u,\psi,x,\phi,y))) &= 
\lambda(u,\psi\phi,x,\phi^\psi,y) \\
&=(x,\phi^\psi,y) \\
(\lambda*1)(u,\psi,x,\phi,y) &= (x,\phi,y)  
\end{align*}
so the condition says that $\phi^\psi=\phi$ whenever $\psi$ has trivial codomain; in other words: 
\begin{itemize}
\item $R_u$ acts on morphisms as the domain functor when $u$ is trivial.
\end{itemize}

\subsection*{The $(\alpha,\rho)$-compatibility condition} This says that the diagram
$$\xymatrix{
X*Y \ar[r]^-{\rho} \ar[dr]_{1*\rho} & (X*Y)*U \ar[d]^{\alpha} \\ &
X*(Y*I) }$$
commutes. Now
\begin{align*}
  \alpha(\rho(x,\phi,y)) &= \alpha(x,\phi,y,1_c,R1_c) \\
  &=(x,\phi 1_c,y,1_c^\phi,y,R1_c) \\
(1*\rho)(x,\phi,y) &= (x,\phi,y,1_c^\phi,R1_c)
\end{align*}
and so the condition says that:
\begin{itemize}
\item the left identity law $\phi 1_{\dom\phi}=\phi$ holds for morphisms $\phi$ in \cc
\item the functors $R_d$ preserve identity morphisms 
\end{itemize}

\subsection*{The $(\lambda,\alpha,\rho)$-compatibility condition}
This says that the composite
$$\xymatrix{
X*Y \ar[r]^-{\rho*1} & (X*U)*Y \ar[r]^{\alpha} & X*(I*Y)
\ar[r]^-{1*\lambda} & X*Y}$$
is the identity. Now
\begin{align*}
(1*\lambda)\alpha(\rho*1)(x,\phi,y) &= (1*\lambda)\alpha(x,1_d,R1_d,\phi,y) \\
&= (1*\lambda)(x,1_d\phi, R1_d,\phi^{1_{\partial(x)}},y) \\
&= (x,1_{d}\phi,y)
\end{align*}
and so the condition says that 
\begin{itemize}
\item the right identity law $1_{\cod\phi} \phi=\phi$ holds for morphisms $\phi$ in \cc.
\end{itemize}

\subsection*{The pentagon.}
This says that the diagram 
$$\xymatrix{
(W*(X*Y))*Z \ar[r]^{\alpha} & W*((X*Y)*Z) \ar[dr]^{1*\alpha} \\
((W*X)*Y)*Z \ar[u]^{\alpha*1} \ar[dr]_{\alpha} && W*(X*(Y*Z)) \\
& (W*X)*(Y*Z) \ar[ur]_{\alpha} }$$
commutes. For $(w,\theta,x,\psi,y,\phi,z)\in ((W*X)*Y)*Z$, we have 
\begin{align*}
(1*\alpha)\alpha  (\alpha*1)(w,\theta,x,\psi,y,\phi,z) &= 
(1*\alpha)\alpha(w,\theta\psi,x,\psi^\theta,y,\phi,z) \\
&= (1*\alpha)(w,(\theta\psi)\phi,x,\psi^\theta,y,\phi^{\theta\psi},z) \\
&= (w,(\theta\psi)\phi,x,\psi^\theta\phi^{\theta\psi},y, (\phi^{\theta\psi})^{\psi^\theta},z) \\
\alpha\alpha(w,\theta,x,\psi,y,\phi,z) &= \alpha(w,\theta,x,\psi\phi,y,\phi^\psi,z) \\
&= (w,\theta(\psi\phi),x,(\psi\phi)^\theta,y,\phi^\psi,z)
\end{align*}
and so the pentagon is equivalent to the conditions
\begin{itemize}
\item $(\theta\psi)\phi = \theta(\psi\phi)$ (associativity of
  composition)
\item $\psi^\theta\phi^{\theta\psi} = (\psi\phi)^\theta$
  (functoriality of $R$)
\item $(\phi^{\theta\psi})^{\psi^\theta} = \phi^\psi$ (double slice condition).
\end{itemize}

I'll summarize these results as:

\begin{theorem}
If \cc is an operadic category and $C$ is its set of objects, there is a skew monoidal category $\sColl_\cc(\Set)$ with underlying category $\Set/C$, with tensor $*$ given by
\[ (X*Y)_c = \sum_{f\colon c\to d} X_d \x \prod_{i\in|d|} Y_{f^{-1}i} ,\]
with unit $U$ consisting of the trivial objects,
 and with structure maps 
\begin{align*}
  \alpha(x,\psi,y,\phi,z) &= (x,\psi\phi,y,\phi^\psi,z) \\
\lambda(u,\phi,x) &= x \\
\rho(x) &= (x,1_{\partial(x)},R1_{\partial(x)})
\end{align*}
\end{theorem}

\begin{example}
In the case of \cp, the object-set is \bbn, and so we obtain a skew
monoidal structure on $\Set/\bbn$ (or equivalently on $\Set^\bbn=[\bbn,\Set]$). By Remarks~\ref{rmk:lambda},~\ref{rmk:rho}, and~\ref{rmk:alpha}, the maps $\alpha$, $\lambda$, and $\rho$ are all invertible, and so we actually have a monoidal structure. This turns out to be use usual substitution monoidal structure on $[\bbn,\Set]$, with respect to which the monoids are the plain operads.
\end{example}

\begin{example}
In the case of \cs, the object-set is once again \bbn, so we obtain another skew monoidal structure on $[\bbn,\Set]$. This time $\lambda$ is invertible, but $\rho$ and $\alpha$ are not, so it is definitely not the same as in the previous example. 
\end{example}

\begin{example}
  If \ca is just a category with object-set $A$, seen as an operadic category with all
  $|a|=\emptyset$ as in Example~\ref{ex:x}, then the corresponding skew monoidal structure on
  $\Set/A$ has unit $0$ and tensor $X*Y$ given by
  $(X*Y)_a=\{(x,\phi)\mid \phi\colon a\to b, x\in
  X_b\}$. Equivalently, this is the restriction of the left Kan extension of $X\colon
  A\to \Set$ along the inclusion $A\to \ca\op$. (This is independent of $Y$.)
\end{example}

\section{Operads} \label{sect:operads}

For any skew monoidal category we can define the category of monoids. In particular we can do so for $\sColl_\cc(\Set)$ for any operadic category \cc. We now unravel what this means.

First we should give an object $T$ of $\sColl_\cc(\Set)$. This amount to giving a set $T_c$ for every object $c\in \cc$.

The unit has the form of a morphism $\eta\colon U\to T$ in $\Set/C$. This amounts to giving, for each trivial object $u\in\cc$, an element $\eta_u\in T_u$.

The multiplication has the form of a morphism $\mu\colon T*T\to T$. This amounts to giving, for each morphism $\phi\colon c\to d$, each $x\in T_d$ and $|d|$-indexed family $y$ with $y_i\in Y_{\phi^{-1}i}$, an element of $T_c$. In other words, for each $\phi\colon c\to d$, we should give a function 
\[ \xymatrix{ T_d \x \prod_{i\in|d|} T_{\phi^{-1}i} \ar[r]^-{\mu(\phi)} &  T_c. } \]

These should satisfy associativity and two unit axioms, which we now explicate. 

Associativity says that the diagram 
\[ \xymatrix{
(T*T)*T \ar[d]_{\alpha} \ar[rr]^{\mu*1} && T*T \ar[d]^{\mu} \\
T*(T*T) \ar[r]_{1*\mu} & T*T \ar[r]_{\mu} & T } \] 
commutes. An element of $(T*T)*T$ involves a pair $(\phi\colon c\to d, \psi\colon d\to e)$; associativity then says that for any such pair, the diagram 
\[ \xymatrix{
T_e \x \prod_{j\in|e|} T_{\psi^{-1}j} \x \prod_{i\in|d|} T_{\phi^{-1}i}  \ar[r]^-{\mu(\phi)\x1} \ar[d]_{\cong} 
& T_d \x \prod_{i\in|d|} T_{\phi^{-1}i} \ar[dd]^{\mu(\psi)} 
\\
T_e\x \prod_{j\in|e|} \left( T_{\psi^{-1}j} \x \prod_{i\in|\psi^{-1}j|} T_{\phi^{-1}i} \right) \ar[d]_-{1\x\prod_j\mu(\phi^\psi_j)} & 
\\
T_e\x \prod_{j\in|e|} T_{(\psi\phi)^{-1}j} 
\ar[r]_{\mu(\psi\phi)} & T_c 
} \] 
diagrams.

The unit condition say that the diagrams
\[ \xymatrix{
U*T \ar[r]^{\eta*1} \ar[dr]_{\lambda} & T*T \ar[d]^{\mu} & T*U \ar[l]_{1*\eta} & T \ar[l]_-{\rho} \ar[dll]^{1} \\
& T 
} \]
commute.

The left unit condition says that for any $\phi\colon c\to u$ with trivial codomain, the diagram 
\[
\xymatrix{
1\x T_c \ar[r]^{\eta\x 1} \ar[dr]_{\cong} & T_u\x T_c \ar[d]^{\mu(\phi)} \\
& T_c } \]
commutes, and the right condition says that for any object $c$ the diagram
\[ 
\xymatrix{
T_c \x \prod_{i\in|c|} T_{1^{-1}_c i} \ar[d]_{\mu(1_c)} & 
T_c \x \prod_{i\in|c|} 1 \ar[l]_{1\x\prod_i \eta} & T_c \ar[l]_-{\cong} \ar[dll]^{1} 
\\
T_c } \]
commutes.

This agrees with the definition in \cite{BataninMarkl-Operadic}, giving:

\begin{theorem}
If \cc is a genuine operadic category then the category $\Op^\cc(\Set)$ of operads defined in \cite{BataninMarkl-Operadic} is the category of monoids in $\sColl_\cc(\Set)$. 
\end{theorem}

The specific structure involved in a \cc-operad for many examples of \cc is spelled out explicitly in \cite{BataninMarkl-Operadic}; in particular, a \cp-operad is a plain operad and an \cs-operad is a symmetric operad. 

\begin{example}
  For a category \ca, seen as an operadic category as in
  Example~\ref{ex:x}, an \ca-operad is precisely a presheaf on \ca. To
  see this, recall that a collection consists of an $A$-indexed family
  $X=(X_a)_{a\in A}$ of sets. The unit of $\sColl_\ca$ is the empty
  set (over $A$), and so any such family $X$ has a unique map $U\to
  X$. The tensor product $X*X$ consists of pairs $(\phi\colon c\to d,
  x\in X_d)$, and in this case $\partial(\phi,x)=c$. Thus to give a
  map $\mu\colon X*X\to X$ is to give, for each $\phi\colon c\to d$ and each $x\in X_d$, an element $x\phi\in X_c$. Associativity of $\mu$ says that $(x\psi)\phi=x(\psi\phi)$, the right unit condition says that $x1_d=x$, and the left unit condition is trivial. 
\end{example}

This gives a cleaner way to regard presheaves as operads than in \cite[Example~1.16]{BataninMarkl-Operadic}, which used $\ca_1$ rather than \ca. An $\ca_1$-operad is a presheaf on \ca, together with an action of a monoid $M$ in \Set.

\section{Functoriality of the {\rm sColl} construction} \label{sect:functoriality}

If $\ce=(\ce,*,U)$ and $\cf=(\cf,*,V)$ are skew monoidal categories, an {\em opmonoidal functor} from \ce to \cf is a functor $F\colon \ce\to \cf$ equipped with a natural transformation $F^2\colon F(X*Y)\to FX*FY$ and a morphism $F^0\colon FU\to V$ satisfying three coherence conditions: one expressing coassociativity of $F^2$ (compatibility with the associativity maps $\alpha$) and two counit conditions for $F^0$ (compatibility with the left unit maps $\lambda$ and the right unit maps $\rho$).
Write $\SkewMon$ for the category of skew monoidal categories and opmonoidal functors.

The goal of this section is to prove the following theorem.

\begin{theorem}\label{thm:functoriality}
There is a functor $\sColl\colon \OpCat\to\SkewMon$ sending an operadic category \cc to the skew monoidal category $\sColl_\cc(\Set)$.
The operadic category \cc is a genuine operadic category in the sense of \cite{BataninMarkl-Operadic} if and only if the left unit constraint $\lambda$ of $\sColl_\cc(\Set)$ is invertible.
\end{theorem}

\proof
Let \cc and \cd be generalized operadic categories. Write $D$ for the set of objects of \cd, and $V$ for the trivial ones, seen as a set over $D$. Let $F\colon \cc\to\cd$ be a strict operadic functor. 

In particular, $F$ determines a function $f\colon C\to D$ between the
sets of objects, and composition with this function induces a functor
$f_!\colon \Set/C\to \Set/D$ which has a right adjoint $f^*$ given by
pullback. For $X\in\Set/C$, the object $f_!(X)$ of $\Set/D$ is just
$X$ as a set, but now with structure map $f\partial$.

Since $f$ maps $U$ to $V$, there is a (unique) map $F^0\colon f_!U\to V$ in $\Set/D$.

Let $X,Y\in\Set/C$. There is a map $F^2\colon f_!(X*Y)\to
f_!(X)*f_!(Y)$ sending $(x,\phi,y)$ to $(x,F\phi,y)$. This is clearly natural in $X$ and $Y$.

We now check the coherence conditions on $F^2$ and $F^0$.
Compatibility with $\lambda$ says that the composite 
$$\xymatrix{
f_!(U*X) \ar[r]^-{F^2} & f_!(U)*f_!(X) \ar[r]^-{F^0*1} & V*f_!(X)
\ar[r]^{\lambda} & f_!(X) }$$
is equal to $f_!(\lambda)$. An element of $f_!(U*X)_a$ has the form
$(u,\phi,x)$, where $\phi\colon c\to u$ is a morphism in \cc with trivial codomain, $x\in X_c$, and  $fc=a$. 
Such an element $(u,\phi,x)$ is sent by $F^2$ to $(u,F\phi,x)$, then
by $F^0*1$ to $(Fu,F\phi,x)$, then by $\lambda$ to $x$; the composite
of these is indeed equal to $f_!(\lambda)$. 

Compatibility with $\rho$ says that the composite 
$$\xymatrix{
f_!(X) \ar[r]^-{f_!(\rho)} & f_!(X*U) \ar[r]^{F^2} & f_!(X)*f_!(U)
\ar[r]^-{1*F^0} & f_!(X)*V}$$
is equal to $\rho$. Now $f_!(\rho)$ sends $x$ to 
$(x,1_{\partial(x)},R1_{\partial(x)})$, then $F^2$ sends this to
$(x,F1_{\partial(x)},R1_{\partial(x)})$, and $1*F^0$ sends this to 
$(x,F1_{\partial(x)},FR1_{\partial(x)})$. Thus compatibility with $\rho$ amounts to the fact that 
\begin{itemize}
\item $F$ preserves identities.
\end{itemize}
(The fact that $FR1_{\partial(x)}=R1_{F\partial(x)}$ then follows.)

Finally, compatibility with $\alpha$ says that the diagram 
$$\xymatrix{
f_!((X*Y)*Z) \ar[r]^{f^2} \ar[d]_{f_!(\alpha)} & 
f_!(X*Y)*f_!(Z) \ar[r]^-{f^2*1} & (f_!(X)*f_!(Y))*f_!(Z)
\ar[d]^{\alpha} \\
f_!(X*(Y*Z)) \ar[r]_{f^2} & f_!(X)*f_!(Y*Z) \ar[r]_-{1*f^2} & f_!(X)*(f_!(Y)*f_!(Z))
}$$
commutes.
An element of $f_!((X*Y)*Z)$ can be written as $(x,\psi,y,\phi,z)$, 
say with $\phi\colon c\to d$, $\psi\colon d\to e$, $x\in X_e$,
and $y$ and $z$ are families indexed by $|e|$ and $|d|$, and
$y_j\in Y_{\psi^{-1}j}$ and $z_i\in Z_{\phi^{-1}i}$. Passing
along the upper path, this gets sent to $(x,\psi,y,F\phi,z)$, then to 
$(x,F\psi,y,F\phi,z)$, then to $(x,F(\psi)F(\phi),y,F(\phi)^{F(\psi)},z)$. 
The lower path goes to $(x,\psi\phi,y,\phi^\psi,z)$, to $(x,F(\psi\phi),y,\phi^\psi,z)$, and
to $(x,F(\psi\phi),y,F(\phi^\psi),x)$. Thus compatibility with $\alpha$ is equivalent to the following two conditions:
\begin{itemize}
\item $F$ preserves composition;
\item $F$ commutes with the $R$ functors. \endproof
\end{itemize}

\begin{remark}\label{rmk:cartesian}
  Observe also that the natural transformation $F^2$ is {\em cartesian}, in the sense that the naturality squares are pullback.
\end{remark}

\section{Collections in other symmetric monoidal categories} \label{sect:V}

In this section we briefly sketch what happens when collections are defined not in \Set but in some other symmetric monoidal category \cv. We suppose that \cv is cocomplete, and that tensoring with a fixed object preserves colimits; this preservation condition certainly holds if \cv is closed.

In the case $\cv=\Set$ it was possible to use $\Set/C$ as the underlying category of collections; for a general \cv, we use $\cv^C$ instead. Apart from this change, everything goes through in essentially the same way. The tensor product is given by the formula
\[
(X* Y)_c = \sum_{\phi\colon c\to d} X_d\ox \bigotimes_{i\in|d|} Y_{\phi^{-1}i} \]
where the sum is over all morphism $\phi\colon c\to d$ in \cc. The unit $U$ is given by 
\[ U_c =
\begin{cases}
  I & \text{if $c$ is trivial} \\
0 & \text{otherwise.}
\end{cases} \]
The $c$-component of the left unit map $\lambda\colon U* X\to X$ has the form 
\[ 
\sum_{\phi\colon c\to d} U_d \ox \bigotimes_{i\in|d|} X_{\phi^{-1}i}\to X_c \]
and is defined via the universal property of the coproduct. If $\phi\colon c\to d$ has trivial codomain then the $\phi$-summand is $I\ox X_c$ and we just use the left unit map $I\ox X_c\to X_c$ in \cv. If $\phi\colon c\to d$ has non-trivial codomain then the $\phi$-summand is $0\ox \bigotimes_i X_{\phi^{-1}i}\cong 0$, and so we use the unique map. The $c$-component of the right unit map $\rho\colon X\to X*U$ has the form 
is given by the composite
\[ \xymatrix{
X_c \ar[r] & X_c\ox \bigotimes\limits_{i\in|c|} I = X_c\ox \bigotimes\limits_{i\in|c|} U_{1^{-1}_c i} \ar[r] & {}\sum\limits_{\phi\colon c\to d} X_d \x \bigotimes\limits_{i\in|d|} U_{\phi^{-1}i} } \]
constructed using right unit maps in \cv and the injection of the $1_c$-summand. The associativity map $\alpha\colon (X*Y)*Z\to X*(Y*Z)$ is more complicated. 
{\allowdisplaybreaks 
\begin{align}
  ((X*Y)*Z)_c &= \sum_{\phi\colon c\to d} (X*Y)_d \ox \bigotimes_{i\in|d|} Z_{\phi^{-1}i}  \notag \\
&= \sum_{\phi\colon c\to d} \left( \sum_{\psi\colon d\to e} X_e\ox \bigotimes_{j\in|e|} Y_{\psi^{-1}j}\right) \ox \bigotimes_{i\in|d|} Z_{\phi^{-1}i} \notag \\
&\cong \sum_{\phi\colon c\to d\atop \psi\colon d\to e} X_e\ox \bigotimes_{j\in|e|} Y_{\psi^{-1}j} \ox \bigotimes_{i\in|d|} Z_{\phi^{-1}i}  \label{eq:left-alpha} \\
&\cong \sum_{\phi\colon c\to d\atop \psi\colon d\to e} X_e\ox \bigotimes_{j\in|e|} \left( Y_{\psi^{-1}j} \ox \bigotimes_{i\in |\psi|^{-1}j} Z_{\phi^{-1}i} \right) \notag \\
(X*(Y*Z))_c &= \sum_{\theta\colon c\to e} X_e\ox \bigotimes_{j\in|e|} (Y*Z)_{\theta^{-1}j} \notag \\
&= \sum_{\theta\colon c\to e} X_e\ox \bigotimes_{j\in|e|} \left(\sum_{\omega_j\colon \theta^{-1}j\to v_j} Y_{v_j}\ox\bigotimes_{i\in|v_j|} Z_{\omega^{-1}_j i}\right) \label{eq:right-alpha}
\end{align}
}
Now $\alpha$ sends the $(\phi,\psi)$-component of \eqref{eq:left-alpha} to the $\psi\phi$-component of \eqref{eq:right-alpha} using the injections
\[ \xymatrix{
Y_{\psi^{-1}j}\ox\bigotimes\limits_{i\in|\psi^{-1}j|} Z_{\phi^{-1}i} \ar[r] & 
{}\sum\limits_{\omega_j\colon \theta^{-1}j\to v_j} Y_{v_j} \ox \bigotimes\limits_{i\in|v_j|} Z_{\omega^{-1}_ji} } \]
of the $\phi^\psi_j$-component for each $j$.

Verification of the five axioms is an exercise in internalizing the arguments already given for the case $\cv=\Set$. This gives a skew monoidal category $\sColl_\cc(\cv)$.

Similarly, if $F\colon \cc\to\cd$ is a strict operadic functor, there is an induced opmonoidal functor $F_!\colon\sColl_\cc(\cv)\to\sColl_\cd(\cv)$, sending the collection $(X_c)_{c\in\cc}$ to $\left(\sum_{Fc=d}X_c\right)_{d\in\cd}$. Unless coproducts in \cv happen to be strictly associative, the construction is now only pseudofunctorial, so that $G_!F_!\cong (GF)_!$, but this does not present any great problem. 

Once again, a \cc-operad in \cv is a monoid in $\sColl_\cc(\cv)$. 

\section{Another characterization of operadic categories} \label{sect:characterization}

We observed above that for any functor $F\colon \cx\to\cc$ and any
object $x\in\cx$, there is an induced functor $F/x\colon \cx/x\to
\cc/Fx$; if \cx has a terminal object $1$, then $\cx/1\cong \cx$, and
we shall allow ourselves to write $F/1$ for the corresponding functor $\cx\to\cc/F1$.
In particular, we could apply this to the functor $\sColl\colon \OpCat\to\SkewMon$, using the fact that \cs is terminal in \OpCat. 

\begin{theorem}\label{thm:coll}
The functor $\sColl/1\colon \OpCat\to\SkewMon/\sColl_{\cs}$ is fully faithful. An opmonoidal functor $P\colon \ce\to\sColl_\cs$ is in the image  of $\sColl/1$ if and only if:
\begin{enumerate}[(a)]
\item the underlying category of \ce is a slice category $\Set/C$,
\item the underlying functor of $P$ is $p_!\colon \Set/C\to\Set/\bbn$ for some $p\colon C\to\bbn$ in \Set, and 
\item the opmonoidal structure map $P^2$ is cartesian. 
\end{enumerate}
\end{theorem}

\proof
Let \cc and \cd be operadic categories, and $P\colon \cc\to\cs$ and $Q\colon \cd\to\cs$ the unique strict operadic functors, and $p\colon C\to\bbn$ and $q\colon D\to \bbn$ their effect on objects. The underlying categories of $\sColl_\cc$ and $\sColl_\cd$ are $\Set/C$ and $\Set/D$, and the underlying functors of $P$ and $Q$ are given by $p_!\colon \Set/C\to\Set/\bbn$ and $q_!\colon \Set/D\to\Set/\bbn$. Any functor $\Set/C\to\Set/D$ commuting with $p_!$ and $q_!$ has the form $f_!$ for a unique function $f\colon C\to D$. 

We are show that to extend the assignment on objects $f\colon C\to D$ to a strict operadic functor $\cc\to\cd$ is equivalent to giving the functor $f_!\colon \Set/C\to\Set/D$ opmonoidal structure which is compatible with those on $p_!$ and $q_!$. The nullary part of the opmonoidal structure amounts to a morphism $f^0\colon f_!(U)\to V$ in $\Set/D$. Such a map is unique if it exists, and will exist if and only if $f$ maps trivial objects to trivial objects; furthermore, when it does exist the compatibility condition is automatic.

The binary part involves natural maps $F^2\colon f_!(X*Y)\to f_!(X)*f_!(Y)$, and the compatibility condition says that the diagram
\begin{equation}\label{eq:f2}
\xymatrix{
q_!f_!(X*Y) \ar@{=}[d] \ar[r]^-{q_!F^2} & q_!(f_!(X)*f_!(Y)) \ar[r]^-{Q^2} & q_!f_!(X)*q_!f_!(Y) \ar@{=}[d] \\
p_!(X*Y) \ar[rr]_{P^2} && p_!(X)*p_!(Y)}  
\end{equation}
should commute. 

We may regard $C$ as an object of $\Set/C$ via the identity morphism $1_C\colon C\to C$. An element of $C*C$ lying over $c$ has the form $(d,\phi,R\phi)$, where $\phi\colon c\to d$ is a morphism in \cc and $R\phi$ is the $|d|$-indexed family consisting of the fibres of $\phi$. The elements of $f_!(C*C)$ are just the same, but now $(d,\phi,R\phi)$ is regarded as lying over $fc$. The elements of $q_!f_!(C*C)=p_!(C*C)$ are still the same, with $(d,\phi,R\phi)$ regarded as lying over $|fc|=|c|$. Finally $P^2\colon p_!(C*C)\to p_!(C)*p_!(C)$ sends $(d,\phi,R\phi)$ to $(d,|\phi|,R\phi)$. Thus a map $F^2\colon f_!(C*C)\to f_!(C)*f_!(C)$ making the relevant instance of the diagram \eqref{eq:f2} commute 
must be of the form $(d,\phi,R\phi)\mapsto (d,F\phi,R\phi)$, for some assignment of a morphism $F\phi\colon fc\to fd$ for every $\phi\colon c\to d$, with $|F\phi|=|\phi|$.

Now let $X$ and $Y$ be arbitrary objects of $\Set/C$. An element of $X*Y$ lying over $c$ has the form $(x,\phi,y)$, where $\phi\colon c\to d$ is a morphism in \cc, $x\in X_d$, and $y$ is a $|d|$-indexed family with $y_i\in Y_{\phi^{-1}i}$. An element of $f_!(X*Y)$ still has the same form, but now $(x,\phi,y)$ is seen as lying over $fc$; and when we pass to $q_!f_!(X*Y)=p_!(X*Y)$ the elements are still unchanged, but now $(x,\phi,y)$ lies over $pc=|c|$. Finally $P^2\colon p_!(X*Y)\to p_!(X)*p_!(Y)$ sends $(x,\phi,y)$ to $(x,|\phi|,y)$.
Thus maps  $F^2\colon f_!(X*Y)\to f_!(X)*f_!(Y)$, natural in $X$ and $Y$ and making the diagram \eqref{eq:f2} 
commute, must have the form $(x,\phi,y)\mapsto (x,F\phi,y)$ for some assignment $(\phi\colon c\to d)\mapsto (F\phi\colon Fc\to Fd)$.

Thus any opmonoidal functor $\sColl_\cc\to\sColl_\cd$ commuting with the induced opmonoidal functors into $\sColl_\cs$ must arise from some assignment $c\mapsto Fc$ and $(\phi\colon c\to d)\mapsto (F\phi\colon Fc\to Fd)$ commuting with the functors $P\colon\cc\to\cs$ and $Q\colon \cd\to\cs$. Just as in the proof of Theorem~\ref{thm:functoriality}, this will define a strict operadic functor $F$ if and only if $F^2$ and $F^0$ satisfy the compatibility conditions for an opmonoidal functor. 

This completes the proof that $\sColl/1$ is fully faithful. It remains to characterize its image.  If $P\colon\ce\to\sColl_\cs$ is in the image then conditions (a) and (b) hold by definition, while (c) holds by Remark~\ref{rmk:cartesian}, $P^2$ is cartesian. Suppose conversely that conditions (a), (b), and (c) all hold. We are given a set $C$ and a function $|~|\colon C\to\bbn$, and we need to try to construct an operadic category \cc with object set $C$.

Define a morphism of \cc to be an element $\phi$ of $C*C$. This will sit over some element $c\in C$ which we define to be its domain. The map $P^2\colon p_!(C*C)\to p_!(C)*p_!(C)$ sends $\phi$ to some element $(d,|\phi|,R\phi)\in p_!(C)*p_!(C)$, where $d\in C$ and $|\phi|\colon |c|\to|d|$ is a function, while $R\phi$ is a $|d|$-indexed family with $|(R\phi)_i|=|\phi|^{-1}i$. Of course we write $\phi^{-1}i$ for the $i$-component of $R\phi$. We define $d$ to be the codomain of $\phi$. Thus we have a directed graph \cc, with a graph morphism into (the underlying directed graph of) \cs sending $c$ to $|c|$ and $\phi$ to $|\phi|$.

Now we use the the fact that $P^2$ is cartesian. For objects $X,Y\in\Set/C$ there are unique maps $X\to C$ and $Y\to C$, and now we have a pullback 
\[ \xymatrix{
p_!(X*Y) \ar[d] \ar[r] & p_!(X)*p_!(Y) \ar[d] \\
p_!(C*C) \ar[r] & p_!(C)*p_!(C) } \] 
in $\Set/\bbn$. An element of $p_!(C*C)$ is still just a morphism $\phi\colon c\to d$, but now it is seen as living over $|c|$ rather than $c$. Its image under $P^2$ is $(d,|\phi|,R\phi)$; thus to give an element of $p_!(X)*p_!(Y)$ living over this is to give $x\in X_d$ and a $|d|$-indexed family $y$ with $y_i\in Y_{\phi^{-1}i}$. This gives the expected description of $X*Y$. To be consistent with the description of $X*Y$, we write $(d,\phi,R\phi)$ for the element of $C*C$ identified with $\phi$. 

An element of $(C*C)*C$ over $c$ has the form $(e,\psi,R\psi,\phi,R\phi)$, where $\phi\colon c\to d$ and $\psi\colon d\to e$. Applying $P^2$ twice sends this to $(e,|\psi|,R\psi,|\phi|,R\phi)$. An element of $C*(C*C)$ over $c$ has the form $(e,\theta,v,\tau,R\tau)$, where $\theta\colon c\to e$ is a morphism in \cc and $\tau\colon R\theta\to v$ is a morphism in $\cc^{|e|}$. Applying $P^2$ twice sends this to $(e,|\theta|,v,|\tau|,R\tau)$. Thus to give $\alpha\colon (C*C)*C\to C*(C*C)$ making the diagram 
\[ \xymatrix{
p_!((C*C)*C) \ar[r]^{P^2} \ar[d]_{p_!(\alpha)} & p_!(C*C)*p_!(C) \ar[r]^{P^2*1} & (p_!C*p_!C)*p_!C \ar[d]^{\alpha} \\
p_!(C*(C*C)) \ar[r]_{P^2} & p_!(C)*p_!(C*C) \ar[r]_{1*P^2} & p_!C*(p_!C*p_!C) } \]
commute is equivalent to giving, for each composable pair $(\phi\colon c\to d, \psi\colon d\to e)$ in \cc, the following data:
\begin{itemize}
\item a morphism $\psi\phi\colon c\to e$ in \cc, with $|\psi\phi|=|\psi||\phi|$;
\item a morphism $\phi^\psi\colon R(\psi\phi)\to R(\psi)$ in $\cc^{|e|}$ with $|\phi^\psi|=|\phi|^{|\psi|}$ and $R(\phi^\psi)=R(\phi)$.
\end{itemize}
By naturality, the map $\alpha\colon (X*Y)*Z\to X*(Y*Z)$ sends $(x,\psi,y,\phi,z)$ to $(x,\psi\phi,y,\phi^\psi,z)$.

Write $U_\cs$ for the unit of $\sColl_\cs$, consisting of the inclusion $\{1\}\to\bbn$.
The unit of $\Set/C$ will be a set $U$ equipped with a morphism $\partial\colon U\to C$. There can be at most one map $p_!U\to U_\cs$, and there will exist one if and only if $|\partial(u)|=1$ for all $u\in U$.

An element of $U*X$ has the form $(u,\phi,x)$, where $u\in U$, $\phi\colon c\to \partial(u)$, and $x\in X_c$. So to give a natural map $\lambda\colon U*X\to X$ whose image under $p_!$ is 
\[ \xymatrix{
p_!(U*X) \ar[r]^-{P^2} & p_!(U)*p_!(X) \ar[r]^-{P^0*1} & V*p_!(X) \ar[r]^-\lambda & p_!(X) } \] 
necessarily has the form $(u,\phi,x)\mapsto x$.

We now seek a natural map $\rho\colon X\to X*U$ for which the composite 
\begin{equation}\label{eq:rho}
 \xymatrix{
p_!X \ar[r]^-{p_!\rho} & p_!(X*U) \ar[r]^-{P^2} & p_!(X)*p_!(U) \ar[r]^-{1*P^0} & p_!(X)*U_\cs } 
\end{equation}
is the $\rho$ for \cs. To give the component at $X=C$ is to give, for each element $c\in C$, a triple $(\cod 1_c,1_c,u(c))$, where $1_c\colon c\to \cod 1_c$ and $u(c)$ is a $|\cod 1_c|$-indexed family of elements of $U$. This will satisfy the condition \eqref{eq:rho} when $\cod 1_c=c$. More generally, if $\rho\colon X\to X*U$ is to be natural and satisfy the condition then it must send $x\in X$ to $(x,1_{\partial(x)},u(\partial(x)))$.

Thus we now have all the data for an operadic category, and we have checked that this induces the skew monoidal structure in the desired way. What remains to be checked are:
\begin{enumerate}[(i)]
\item $\partial U\to C$ is injective
\item associativity laws for  \cc
\item identity laws for \cc
\item functoriality of $R$
\item $R_{\partial(u)}=\dom$ if $u\in U$
\item the double slice condition for $R$.
\end{enumerate}

All but the first of these follow as in Section~\ref{sect:sColl}, when
we checked that axioms for the skew monoidal category $\sColl_\cc$,
since in each case we proved the equivalence of the axiom with some
subset of the conditions above. As for (i),  we may use the
$(\lambda,\rho)$-compatibility condition which says that the composite 
\[ \xymatrix{ U \ar[r]^-{\rho} & U*U \ar[r]^-{\lambda} & U } \]
is the identity. But for $v\in U$ we have $\rho(v)=(\partial(v),1_{\partial(v)},u(\partial(v)))$, and $\lambda$ sends this to $u(\partial(v))$. Now $u(\partial(v))$ depends only on $\partial(v)$; thus if it is to be equal to $v$, then $\partial$ must indeed be injective, giving (i).
\endproof

\begin{remark}
Let \cc be an operadic category for which $|c|=1$ for all $c\in \cc$. Then the operadic structure consists of a functor $R\colon \cc/c\to \cc$ for each $c\in \cc$.
The functor $\Set/C\x\Set/C\to\Set/C$ preserves connected limits, and so corresponds to a {\em span} $m$ from $C\x C$ to $C$ (with vertex $C*C$). Similarly the unit $U\to C$ can be seen as a span $i$ from $1$ to $C$. Furthermore, the structure maps $\alpha$, $\lambda$, and $\rho$ can be seen as morphisms of  spans
\[ \xymatrix{
C\x C\x C \ar[r]^{m\x 1}_{~}="1" \ar[d]_{1\x m} & C\x C \ar[d]^m \\
C\x C \ar[r]_{m}^{~}="2" & C \ar@{=>}"1";"2"^{\alpha} }\quad
\xymatrix{
1\x C \ar[r]^{i\x 1}_{~}="3" \ar@{=}[d] & C\x C \ar[d]^{m} \\
1\x C \ar[r]_{\pi_2}^{~}="4" & C 
\ar@{=>}"3";"4"^{\lambda}
}\quad
\xymatrix{
C\x C \ar[d]^{m} & C\x 1 \ar[l]_{1\x i}^{~}="5" \ar@{=}[d] \\
C & C\x 1 \ar[l]^{\pi_1}_{~}="6"
\ar@{=>}"6";"5"_{\rho} }
 \] 
and the axioms say that these define a skew monoidale (internal skew
monoidal structure) in the monoidal bicategory \Span. These
observations allow us to recover from Theorem~\ref{thm:coll} the
characterization of skew monoidales in \Span given in
\cite[Chapter~3]{Andrianopoulos-thesis}. 

It is also possible to characterize $\sColl_\cc$ in terms of \Span for
a general operadic category \cc; this involves using \cs to define a
skew monoidal bicategory structure on $\Span^\bbn$, then considering
skew monoidales in this skew monoidal bicategory. 
\end{remark}

\section{Fibrewise trivial morphisms} \label{sect:bbc}

Let \cc be a generalized operadic category. Recall that $\phi\colon c\to d$ is said to be 
{\em fibrewise trivial} when all fibres of $\phi$ are in $U$. 
By assumption, every identity morphism is fibrewise trivial; we shall see in Lemma~\ref{lemma:triviality-reflection} that the fibrewise trivial morphisms are closed under composition, and so there is a category \bbc consisting of all objects of \cc and fibrewise morphisms between them. By extension, we define a morphism $\pi$ in $\cc/e$ to be fibrewise trivial if $\dom(\pi)$ is fibrewise trivial in \cc; and we define a morphism $\pi$ in $\cc^{|e|}$ to be fibrewise trivial if $\pi_i$ is fibrewise trivial in \cc for all $i\in |e|$.

Since $|\pi^{-1}i|=|\pi|^{-1}i$, if $\pi$ is fibrewise trivial then so is $|\pi|$.

\begin{lemma}\label{lemma:Rbbc}
  $R\colon \cc/e\to \cc^{|e|}$ preserves fibrewise triviality. 
\end{lemma}

\proof
Let $\phi\colon (\psi\phi\colon c\to e)\to (\psi\colon d\to e)$ be fibrewise trivial in $\cc/e$, so that $\phi\colon c\to d$ is fibrewise trivial in \cc. We are to show that $\phi^\psi_j\colon (\psi\phi)^{-1}j\to \psi^{-1}j$ is fibrewise trivial for each $j\in |e|$. But each fibre $(\phi^\psi_j)^{-1}i$ of $\phi^\psi_j$ is just a fibre $\phi^{-1}i$ of $\phi$, thus trivial. \endproof

\begin{lemma}\label{lemma:triviality-reflection}
If $\pi\colon d\to e$ is fibrewise trivial, then a morphism
$\phi\colon c\to d$ is fibrewise trivial if and only if $\pi\phi$ is so.  
\end{lemma}

\proof
For each $j\in |e|$ we have a morphism
$$\xymatrix{
(\pi\phi)^{-1}j \ar[r]^{\phi^\pi_j} & \pi^{-1}j }$$
and $(\phi^\pi_j)^{-1}i=\phi^{-1}i$ for all $i\in |\pi^{-1}j|$.

Since $\pi$ is fibrewise trivial, there is a unique $i\in|\pi^{-1}j|$, and since the codomain $\pi^{-1}j$ of $\phi^\pi_j$ is trivial, the fibres of $\phi^\pi_j$ are trivial.

Combining these two facts, we see that
$(\pi\phi)^{-1}(|\pi|i)=\phi^{-1}i$ for all $i\in|d|$, and so that
$\pi\phi$ has trivial fibres if and only if $\phi$ does so.
\endproof

In particular, the fibrewise trivial morphisms form a subcategory of
\cc, which I'll call \bbc. Furthermore, the inclusion reflects
isomorphisms: any isomorphism in \cc which is fibrewise trivial also
has fibrewise trivial inverse.
On the other hand, the inclusion need not be full on isomorphisms: there can be isomorphisms in \cc which are not fibrewise trivial. 





\begin{example}
If \cc is the operadic category \cp for plain operads, then the only fibrewise trivial morphisms are the identities and so the category $\bbc$ is still just the discrete category $\bbn$. 
\end{example}

\begin{example}
If \cc is the operadic category \cs for symmetric operads, then the fibrewise trivial morphisms are the bijections, and so $\bbc$ is the category \bbp of finite sets and bijections. 
\end{example}

\begin{example}
If \cc is a category \ca, seen as an operadic category with $|a|=0$ for all $a$, then all morphisms are fibrewise trivial, and so $\bbc$ is just \ca once again. 
\end{example}

\section{Normalization} \label{sect:normalization}

If $\ce=(\ce,*,U)$  is a skew monoidal category which has coequalizers of reflexive pairs, and tensoring on the right with any object preserves these coequalizers, then there is a way \cite{mw} to associate a right-normal skew monoidal category $\ce^U$ with the same category of monoids as \ce, which we now recall. 

The unit object $U$ of a skew monoidal category has a monoid structure with multiplication $\lambda\colon U*U\to U$ and unit $1\colon U\to U$. We may define a right $U$-module to be an object $X\in\ce$ equipped with a map $r_X\colon X*U\to X$, or just $r$, called the action and  making the diagrams 
\[
\xymatrix{
(X*U)*U \ar[r]^{\alpha} \ar[d]_{r*1} & X*(U*U) \ar[r]^{1*\lambda} & X*U \ar[d]^{r} & X \ar[r]^{\rho} \ar[dr]_{1}  & X*U \ar[d]^{r} \\
X*U \ar[rr]_{r} && X && X } \]
commute. A homomorphism of $U$-modules is a morphism between the underlying objects which commutes in the obvious sense with the actions. These $U$-modules and their homomorphisms form the category $\cc^U$. 

In fact every object $X$ of \cc has a canonical {\em left} action $\lambda\colon U*X\to X$, preserved by any morphism, and ``compatible'' with any right action, so that $\cc^U$ can also be thought of as the category of {\em bimodules} over $U$. This suggests that there should be a tensor product on $\cc^U$ defined as for bimodules using a coequalizer, and this is indeed the case. Explicitly, the tensor product of $U$-modules $(X,r_X)$ and $(Y,r_Y)$ is given by the coequalizer 
\[ \xymatrix @R1pc {
(X*U)*Y \ar[rr]^{r_X*1} \ar[dr]_{\alpha} && X*Y \ar[r]^{p} &  X\wedge Y \\
& X*(U*Y) \ar[ur]_{1*\lambda} } \]
in \cc, and this becomes a $U$-module via a map $r\colon (X\wedge Y)*U\to X\wedge Y$ which is uniquely determined by commutativity of 
\[ \xymatrix{
(X*Y)*U \ar[r]^{p*1} \ar[d]_{\alpha} & (X\wedge Y)*U \ar@{.>}[d]^{r} \\
X*(Y*U) \ar[r]_{1*r_Y} & X*Y. } \]
The forgetful functor $\cc^U\to \cc$ has an opmonoidal structure defined using the maps $p\colon X*Y\to X\wedge Y$. See \cite {mw} for further details. 

We shall apply this in the case $\ce=\sColl_\cc$ for an operadic category \cc. As a category $\sColl_\cc$ is just $\Set/C$ which is of course cocomplete; and tensoring on the right with an object is cocontinuous, so we can form $\sColl^U_\cc$. 

\begin{proposition}
The category $\sColl^U_\cc$ is equivalent to the category of presheaves on \bbc.
\end{proposition}

\proof
We use the fact that $\Set/C\simeq [C,\Set]$, and show that a $U$-action on $X\to C$ turns the corresponding $X\colon C\to \Set\colon c\mapsto X_c$ into a presheaf. 

Let $X$ be an object of $\Set/C$. An element of $(X*U)_c$ is a triple $(x,\phi,u)$, where $\phi\colon c\to d$ is a morphism in \cc, $x\in X_d$, and $u$ is a $|d|$-indexed family with $u_j\in U_{\phi^{-1}j}$. But there is at most one element in $U_{\phi^{-1}j}$, and there will be one if and only if $\phi^{-1}j$ is trivial. Thus in fact an element of $(X*U)_c$ amounts to a fibrewise trivial morphism $\phi\colon c\to d$ and an element $x\in X_d$. 

Thus for any $X$ in $\Set/C$, to give a morphism $X*U\to X$ is to give, for each $\pi\colon c\to d$ in \bbc and each $x\in X_d$, an element $x\pi\in X_c$.  The associativity and unit conditions for a $U$-module say that $(x\pi)\sigma=x(\pi\sigma)$ and $x1_d=x$. Thus a $U$-module is precisely a presheaf on \bbc; furthermore, the condition for a map $X\to Y$ to preserve the action is precisely naturality. 
\endproof

\begin{proposition}
  $U$ has a unique $U$-module structure, and this is the unit for $\sColl^U_\cc$. The tensor product $X\wedge Y$ of $U$-modules $X$ and $Y$ is the quotient of $X*Y$ by the equivalence relation generated by 
\[ (x,\pi\phi,y) \sim (x\pi,\phi,y_\pi) \]
where $\phi\colon c\to d'$ is in \cc and $\pi\colon d'\to d$ is in \bbc, where $x\in X_d$ and $y=(y_i\in Y_{(\pi\phi)^{-1}i})_{i\in|d|}$, and where $y_\pi$ is the $|d'|$-indexed family with $(y_\pi)_{i'}=y_{\pi i'}$. If $\sigma\colon c'\to c$ is in \bbc, then 
\[ [x,\phi,y]\sigma = [x,\phi\sigma,y\sigma]\]
where $\phi\colon c\to d$, $x\in X_d$, and $y=(y_i\in Y_{\phi^{-1}i})_{i\in |d|}$; and where $y\sigma$ is the $|d|$-indexed family with $(y\sigma)_i=y_i \sigma^\phi_i$.
\end{proposition}

\proof
This is essentially all true by construction. We just point out that $[x,\phi,y]$ denotes the equivalence class of $(x,\phi,y)$, and that in the last part, $\sigma^\phi_i$ is fibrewise trivial because $\sigma$ is so, thus the presheaf structure of $Y$ allows us to form $y_i\sigma^\phi_i$. 
\endproof

\begin{example}
If \cc is the operadic category \cp for plain operads, then $\bbc=C=\bbn$; in this case $\sColl_\cp$ is the usual monoidal category $[\bbn,\Set]$ of plain collections, and in particular is already right normal; thus $\sColl^U_\cp=\sColl_\cp$.
\end{example}

\begin{example}
If \cc is the operadic category \cs for symmetric operads, then $\bbc$ is the category \bbp of finite sets and bijections. Furthermore $\sColl^U_\cs$ is the usual monoidal category of collections $[\bbp,\Set]$.
\end{example}

\begin{example}
If \cc is a category \ca, seen as an operadic category with $|a|=0$ for all $a$, then all morphisms are fibrewise trivial, and so $\bbc$ is just \ca once again. An element of $X*Y$ consists of a $\phi\colon c\to d$ in \ca, and an element $x\in X_d$; we write such an element as $(x,\phi)$. An element of $X\wedge Y$ (for any $Y$) is an equivalence class $[x,\phi]$ of the equivalence relation on $X*Y$ generated by $(x,\psi\phi)\sim (x\psi,\phi)$. Clearly $(x,\phi)\sim(x',\phi')$ if and only if $x\phi= x'\phi'$, thus each $X\wedge Y$ is canonically isomorphic to $X$ itself. 
\end{example}

The general theory guarantees that $\sColl^U_\cc$ is right normal; we now investigate when it is left normal and when it is Hopf. For the first of these there is an easy necessary and sufficient condition. Recall that a functor $F\colon \ca\to\cb$ and an object $B\in\cb$, there is a category $B/F$ whose objects are objects $A\in\ca$ equipped with a morphism $\phi\colon B\to FA$ in \cb, and whose morphisms $(A,\phi)\to (A',\phi')$ are morphisms $\psi\colon A\to A'$ in \ca for which $F\psi.\phi=\phi'$. Then $F$ is said to be {\em final} if $B/F$ is connected for all $B\in\cb$.

\begin{proposition}
  For an operadic category \cc, the skew monoidal category $\sColl^U_\cc$ is left normal if and only if the inclusion $\bbc\to\cc$ is final. This will always be the case if \cc is a genuine operadic category in the sense of \cite{BataninMarkl-Operadic}.
\end{proposition}

\proof
It is true in general that if \ce is a left normal skew monoidal category satisfying the conditions for $\ce^U$ to exist, then $\ce^U$ will also be left normal \cite{mw}. This implies the second sentence, but we shall also see this on our way to proving the first.

An element of $(U*X)_c$ consists of a morphism $\phi\colon c\to u$ in
\cc with $u$ trivial, and an element $x\in X_c$. Then $\lambda$ maps
the pair $(\phi,x)$ to $x$. Of course this is invertible (for all $X$)
if and only if, for each $c\in C$ there is a unique map to some
trivial object (in other words \cc is a genuine operadic category).

On the other hand $U\wedge X\to X$ will be invertible for all $X$ when
\begin{itemize}
\item for each $c\in\cc$ there exists a morphism $c\to u$ with $u$
  trivial
\item for each $c\in \cc$, any two morphisms to a trivial object are
  equivalent under the equivalence relation $\phi\sim\pi\phi$, where
  $\pi$ is a morphism in \bbc between trivial objects.
\end{itemize}
But this says precisely that the comma category is connected.
\endproof

Next we turn to the Hopf condition (invertibility of $\alpha$). In this case we give a sufficient condition only. This condition is a bit more complicated to state, and so to motivate it we start with the following definition of weak right adjoint. The word ``weak'' refers to a universal property which has been weakened to involve only existence rather than uniqueness. There is more than one notion which might reasonably be given this name, but this is the one which is most useful in this paper.

\begin{definition}
  Let $R\colon \ca\to\cb$ be a functor. A{\em weak right adjoint to $R$} to consists of the
  following:
  \begin{enumerate}[(i)]
  \item for each object $B\in \cb$ an object $SB\in\ca$ and a morphism
    $\sigma_B\colon RSB\to B$ in \cb
\item for each object $A\in\ca$ and each $\tau\colon RA\to B$ in \cb, a
  morphism $\pi\colon A\to SB$ for which the composite
\[\xymatrix{ RA \ar[r]^{R\pi} & RSB \ar[r]^-{\sigma_B} & B }\]
is equal to $\tau$.
  \end{enumerate}
The $\sigma_B$ will be called the (components of the) counit. In
particular, for every $A\in\ca$ there is a morphism $\pi_A\colon A\to
SRA$ making the diagram 
$$\xymatrix @R1pc {
RA \ar[dr]^{1} \ar[dd]_{R\pi_A} \\
& RA \\
RSRA \ar[ur]_{\sigma_{RA}} }$$
commute. This might be called ``the unit'', although it is not
uniquely determined by the remaining data. 
\end{definition}

We have seen that a functor $F\colon \ca\to\cb$ induces a functor $F/a\colon \ca/a\to \cb/Fa$ for any object $a\in\ca$. Similarly it induces a functor $a/F\colon a/\ca\to Fa/\cb$.

\begin{proposition}\label{prop:Hopf}
Let \cc be an operadic category for which $\theta/R$ has a
weak 
right adjoint $S_\theta$ for each $\theta\colon c\to e$, and for which
the components of the counits of these weak right adjoints are fibrewise
trivial. Then the induced skew monoidal structure on $[\bbc\op,\Set]$
is Hopf (it has an invertible associativity map).  
\end{proposition}

\proof
An element of $(X*Y)*Z$ has the form $(x,\psi,y,\phi,z)$ where
$\phi\colon c\to d$, $\psi\colon d\to e$, $x\in X_e$, $y=(y_j\in
Y_{\psi^{-1}j})_{j\in|e|}$, $z=(z_i\in
Z_{\phi^{-1}i})_{i\in|d|}$. This has an action of $\bbc$ where for
$\pi\colon c'\to c$
$$(x,\psi,y,\phi,z)\pi = (x,\psi,y,\phi\pi,z(\pi^\phi)).$$
We obtain $(X\wedge Y)\wedge Z$ from $(X*Y)*Z$ by factoring out by the
equivalence relation generated by 
$$(x,\sigma\psi,y,\pi\phi,z) \sim
(x\sigma,\psi\pi,y_\sigma(\pi^\psi),\phi,z_\pi)$$
where 
$$\xymatrix{
c \ar[r]^{\phi} & d' \ar[r]^{\pi} & d \ar[r]^{\psi} & e'
\ar[r]^{\sigma} & e }$$
with $\pi$ and $\sigma$ in \bbc; where $x\in X_{e}$, $y=(y_j\in
Y_{(\sigma\psi)^{-1}j})_{j\in|e|}$, and where $z=(z_i\in
Z_{(\pi\phi)^{-1}i})_{i\in |d|}$.

An element of $X*(Y*Z)$ has the form $(x,\theta,(y,\omega,z))$, where
$\theta\colon c\to f$ and $x\in X_f$; where $y=(y_j\in
Y_{\theta^{-1}j})_{j\in|f|}$ and $\omega=(\omega_j\colon
\theta^{-1}j\to v_j)_{j\in|e|}$, and where $z=(z_{i,j}\in
Z_{\omega^{-1}_ji})_{j\in|e|, i\in|v_j|}$. We obtain $X\wedge(Y\wedge
Z)$ by factoring out by the equivalence relation generated by 
$$(x,\sigma\theta,(y,\tau\omega,z)) \sim
(x\sigma,\theta,(y_\sigma\tau_\sigma,\omega_\sigma,z_{\sigma,\tau_\sigma}))$$
where 
$$\xymatrix{
c \ar[r]^{\theta} & f' \ar[r]^{\sigma} & f & \theta^{-1}j
\ar[r]^{\omega_j} & v'_j \ar[r]^{\tau_j} & v_j }$$
with $\sigma$ and the $\tau_j$ in \bbc; where $x\in X_{f}$, 
$y=(y_j\in Y_{(\sigma\theta)^{-1}j})_{j\in|f|}$, and where 
$z=(z_{j,i}\in Z_{(\tau_j\omega_j)^{-1}i})_{j\in|f|,i\in|v_j|}$. Here
$(y_{\sigma})_{j'}=y_{|\sigma|j'}$, 
$(\tau_{\sigma})_{j'}=\tau_{|\sigma|j'}$, and
$(z_{\sigma,\tau_\sigma})_{j',i'} = z_{|\sigma|j',|\tau_{|\sigma|j'}|i'}$.

Of course $\alpha$ sends the equivalence class $[x,\psi,y,\phi,z]$ to 
$[x,\psi\phi,[y,\phi^\psi,z]]$.

Let $\theta\colon c\to e$ in $\cc$ and $\omega\colon R(\theta)\to v$
in $\cc^{|e|}$ be given. By assumption, the functor $\theta/R\colon \theta/\cc/e\to R\theta/\cc^{|e|}$ has a weak right adjoint $S_\theta$. Applying this to $\omega$ gives a
factorization
$$\xymatrix{ c \ar[r]^{\phi} & d \ar[r]^{\psi} & e}$$
of $\theta$. The counit gives a factorization 
$$\xymatrix{
R\theta \ar[r]^{\phi^{\psi}} \ar[dr]_{\omega} & R\psi
\ar@{.>}[d]^{\sigma} \\ & v }$$
or, in terms of components,
$$\xymatrix{
(\psi\phi)^{-1}j \ar[dr]_{\omega_j} \ar[r]^{\phi^\psi_j} & \psi^{-1}j
\ar@{.>}[d]^{\sigma_j} \\ & v_j }$$
for each $j$, with these $\sigma_j$ lying in \bbc. 

The (weak) universal property guarantees that for any other factorization
$\theta=\psi'\phi'$, equipped with a morphism $\tau$ as in 
$$\xymatrix @R1pc {
& R(\psi') \ar[dd]^{\tau} \\
R(\theta) \ar[ur]^{(\phi')^{(\psi')}} \ar[dr]_{\omega} \\
& v }$$
there is a morphism $\pi$ for which the diagrams 
$$\xymatrix @R1pc {
& d' \ar[dr]^{\psi'} \ar[dd]^{\pi} & & R(\psi') \ar[dd]_{\pi^\psi}
\ar[dr]^{\tau} \\
c \ar[ur]^{\phi'} \ar[dr]_{\phi} && e & & v\\
& d \ar[ur]_{\psi} && R(\psi) \ar[ur]_{\sigma} }$$
commute. 

In fact we only need $\pi$ to exist in the case where the components of $\tau$ are
fibrewise-trivial. We do need, however, that in that case the
components of $\pi$ will also be fibrewise-trivial. 

To see that that is true, observe that since $\sigma$ is fibrewise-trivial, we know by
Lemma~\ref{lemma:triviality-reflection} that $\tau$ is
fibrewise-trivial if and only if $\pi^\psi$ is so. But
$(\pi^\psi)^{-1}i=\pi^{-1}i$, and so $\pi^\psi$ is fibrewise-trivial
if and only if $\pi$ is so. In other words, the adjointness respects
fibrewise-triviality. 

Now consider  $(x,\theta,(y,\omega,z))\in X*(Y*Z)$, with $\theta\colon
c\to e$ and $\omega\colon R\theta\to v$. 
For each $j$ we have $y_j\in v_j$, and so a $y_j\sigma_j\in
Y_{\psi^{-1}j}$. Writing $y\sigma$ for the family
$(y_j\sigma_j)_{j\in|f|}$, we get an object $(x,\psi,y\sigma)\in X*Y$.

For each $j\in |f|$ and each $i\in|v_j|$ we have a $z_i\in
Z_{\omega^{-1}_ji}$.
If now $i\in|\psi^{-1}j|$, we have $\phi^{-1}i=(\phi^\psi)^{-1}_j i =
\omega^{-1}_j(|\sigma_j|i)$, and so we obtain a $|d|$-indexed family
$z_\sigma$ with $(z_{\sigma})i = z_{j,|\sigma_j|i}$, where $j=|\psi|i$.
Thus in fact we have
$(x,\psi,y\sigma,\phi,z_{\sigma})\in (X*Y)*Z$. This defines a function 
$$\xymatrix @R0pc {
X*(Y*Z) \ar[r]^{\alpha'} & (X*Y)*Z \\
(x,\theta,(y,\omega,z)) \ar@{|->}[r] &
(x,\psi,y\sigma,\phi,z_{\sigma}).
}$$
Now
\begin{align*}
  \alpha(\alpha'(x,\theta,(y,\omega,z)) &=
                                          \alpha(x,\psi,y\sigma,\phi,z_\sigma)
  \\
&= (x,\psi\phi,(y\sigma,\phi^\psi,z_\sigma)) \\
&\sim (x,\psi\phi,(y,\sigma\phi^\psi,z)) \\
&= (x,\psi\phi,(y,\omega,z))
\end{align*}
and so 
$\alpha\colon (X\wedge Y)\wedge Z\to X\wedge(Y\wedge Z)$ is surjective.

We know that $\alpha$ passes to the quotient in full generality, but
we need to see what happens to $\alpha'$. 
Suppose then that we have two elements
$(x,\theta,(y,\tau\omega,z))$ and
$(x,\theta,(y\tau,\omega,z_{1,\tau})$
of $X*(Y*Z)$ which are related by $\tau$.
Applying  the right adjoint to 
$$\xymatrix{
& v' \ar[d]^{\tau} \\
R\theta \ar[ur]^{\omega} \ar[r]_{\tau\omega} & v }$$
gives
$$\xymatrix @R1pc {
& d' \ar[dr]^{\psi'} \ar[dd]^{\pi} \\
c \ar[ur]^{\phi'} \ar[dr]_{\phi} && f \\
& d \ar[ur]_{\psi} }$$
and now we have 
$$\xymatrix @R1pc {
& R(\psi') \ar[r]_{\sigma'} \ar[dd]^{\pi^\psi} & v' \ar[dd]^{\tau} \\
R(\theta) \ar[ur]_{(\phi')^{(\psi')}} \ar[dr]^{\phi^\psi}
\ar@/^3pc/[urr]^{\omega} \ar@/_3pc/[drr]_{\tau\omega} \\
& R(\psi) \ar[r]^{\sigma} & v }$$

Now
\begin{align*}
\alpha'(x,\theta,(y,\tau\omega,z)) &= (x,\psi,y\sigma,\phi,z_\sigma) \\
  &= (x,\psi,y\sigma,\pi\phi',z_\sigma) \\
&\sim (x,\psi\pi,y\sigma\pi^\psi,\phi',z_{\sigma\pi^\psi}) \\
&=(x,\psi',y\tau\sigma',\phi',z_{\tau\sigma'}) \\ 
&=\alpha'(x,\theta,(y\tau,\omega,z_{(1,\tau)}))
\end{align*}
as desired.
Next we compare $\alpha'(x,\tau\theta,(y,\omega,z))$ and
$\alpha'(x\tau,\theta,(y_\tau,\omega_\tau,z_{1,\tau}))$.

Applying the right adjoint for $\tau\theta$ to $\omega\colon
R(\tau\theta)\to v$ gives a factorization $(\phi,\psi)$, while
applying the right adjoint for $\theta$ to $\omega_\tau\colon
R(\theta)=R(\tau\theta)_\tau\to v_\tau$ gives a factorization
$(\phi',\psi')$, as in 
$$\xymatrix @R1pc {
& d' \ar[r]^{\psi'} & f' \ar[dd]^{\tau} \\
c \ar[ur]^{\phi'} \ar[dr]_{\phi} \ar[urr]_{\theta} \\
& d \ar[r]_{\psi} & f }$$
%
and now the counits give
$$\xymatrix @R1pc {
& R\psi' \ar[dd]^{\sigma'} & & R\psi \ar[dd]^{\sigma} \\
R\theta \ar[ur]^{(\phi')^{\psi'}} \ar[dr]_{\omega_\tau} && R(\tau\theta) \ar[ur]^{\phi^\psi} \ar[dr]_{\omega} \\
&v_\tau && v .}$$
The first of these can equally be seen as having the form
$$\xymatrix @R1pc {
& R(\tau\psi')_\tau \ar[dd]^{\sigma'} \\
R(\tau\theta)_\tau \ar[ur]^{(\phi')^{(\tau\psi')}_\tau}
\ar[dr]_{\omega_\tau} \\
& v_\tau }$$
and so determines a unique
$$\xymatrix @R1pc {
& R(\tau\psi') \ar[dd]^{\sigma''} \\
R(\tau\theta) \ar[ur]^{(\phi')^{(\tau\psi')}} \ar[dr]_{\omega} \\
& v }$$
where $\sigma''_\tau=\sigma'$, 
and now by the ``universal property of $(\phi,\psi)$'' there is a
fibrewise-trivial $\pi$ making the diagrams 
$$\xymatrix @R1pc {
& d' \ar[dr]^{\tau\psi'} \ar@{.>}[dd]^{\pi} \\
c \ar[ur]^{\phi'} \ar[dr]_{\phi} && f \\
& d \ar[ur]_{\psi} }\quad
\xymatrix @R1pc {
R(\tau\psi') \ar[r]^{\pi^\psi} \ar[ddr]_{\sigma''} & R(\psi)
\ar[dd]^{\sigma} \\ \\ & v }$$
commute, and now
\begin{align*}
  \alpha'(x,\tau\theta,(y,\omega,z)) &= (x,\psi,y\sigma,\phi,z_\sigma)
  \\
&= (x,\psi,y\sigma,\pi\phi',z_\sigma) \\
&\sim (x,\psi\pi,y\sigma\pi^\psi,\phi',z_{\sigma\pi^\psi}) \\
&= (x,\tau\psi',y\sigma'',\phi',z_{\sigma''} \\
&\sim (x\tau,\psi',y_\tau\sigma',\phi',z_{\sigma''}) \\
&= \alpha'(x\tau,\theta,(y_\tau,\omega_\tau,z_{1,\tau}))
\end{align*}
and so $\alpha'$ is indeed well-defined as a map $X\wedge(Y\wedge
Z)\to (X\wedge Y)\wedge Z$. 

Finally, it remains to show that at the level of the quotients,
$\alpha'\alpha$ is the identity. Given $\phi\colon c\to d$ and
$\psi\colon d\to e$, we may form $\psi\phi\colon c\to e$ and
$\phi^\psi\colon R(\psi\phi)\to R(\psi)$, and now the right adjoint
gives another factorization $(\phi',\psi')$ and the counit has the form
$$\xymatrix @R1pc {
& R(\psi') \ar[dd]^{\sigma} \\
R(\psi\phi) \ar[ur]^{(\phi')^{(\psi')}} \ar[dr]_{\phi^\psi}  \\
& R(\psi) }$$
and the unit $\pi$ has the form
$$\xymatrix{
& d \ar[dr]^{\psi} \ar[dd]^{\pi} \\
c \ar[ur]^{\phi} \ar[dr]_{\phi'} && e \\
& d' \ar[ur]_{\psi'} }$$
and one of the triangle equations says that the composite
$$\xymatrix{
R(\psi) = R(\psi'\pi) \ar[r]^-{\pi^{\psi'}} & R(\psi') \ar[r]^{\sigma}
& R(\psi) }$$
is the identity. Now
\begin{align*}
  \alpha'(\alpha(x,\psi,y,\phi,z)) &= \alpha'(
                                     x,\psi\phi,(y,\phi^\psi,z)) \\
&= (x,\psi',y\sigma,\phi',z_\sigma) \\
&= (x,\psi',y\sigma,\pi\phi,z_\sigma) \\
&\sim (x,\psi'\pi,y\sigma\pi^{\psi'},\phi,z_{\sigma\pi^{\psi'}}) \\
&= (x,\psi,y,\phi,z)
\end{align*}
as required. 
\endproof

As observed in the proof of the proposition, we do not actually need the full strength of the universal property of a weak right adjoint. We record this observation as the following proposition.

\begin{proposition}\label{prop:Hopf2}
  Let \cc be an operadic category. The induced skew monoidal structure on $[\bbc\op,\Set]$ will be Hopf provided that the following conditions are satisfied:
  \begin{itemize}
  \item for each $\theta\colon c\to e$ and each fibrewise trivial $\omega\colon R(\theta)\to v$ there is a factorization 
\[ \xymatrix{ c \ar[r]^{\phi} & d\ar[r]^{\psi} & e } \]
of $\theta$, and a factorization 
\[ \xymatrix{ R\theta \ar[r]^{\phi^\psi} & R\psi \ar[r]^{\sigma} & v } \]
of $\omega$, with $\sigma$ fibrewise trivial;
\item if 
\[ \xymatrix @R1pc { c \ar[r]^{\phi'} & d' \ar[r]^{\psi'} & e \\
R(\theta) \ar[r]^{(\phi')^{(\psi')}} & R(\psi') \ar[r]^{\tau} & v } \]
are also factorizations of $\theta$ and $\omega$, then there is a morphism $\pi\colon d'\to d$ making the diagrams
$$\xymatrix @R1pc {
& d' \ar[dr]^{\psi'} \ar[dd]^{\pi} & & R(\psi') \ar[dd]_{\pi^\psi}
\ar[dr]^{\tau} \\
c \ar[ur]^{\phi'} \ar[dr]_{\phi} && e & & v\\
& d \ar[ur]_{\psi} && R(\psi) \ar[ur]_{\sigma} }$$
commute. 
  \end{itemize}
\end{proposition}

\section{Examples} \label{sect:examples}

We have already discussed the basic examples \cs (for symmetric operads) and \cp (for plain operads). In the case of \cp, the skew monoidal category $\sColl_\cp$ is already monoidal, and in fact the $R$ functors are isomorphisms, thus so too are the $\theta/R$.
In the case of \cs, the $R$ functors are surjective equivalences, thus so too are the $\theta/R$, and so $\sColl^U_\cs$ is monoidal; indeed it is equivalent to the monoidal category $\Coll(\Set)=[\bbp,\Set]$.

In the case of a category \ca regarded as an operadic category with $|a|=0$ for all $a\in\ca$, the $R$-functors have the form $\ca/a\to 1$. These have right adjoints with identity counit, thus the same is true of the $\theta/R$, and so $\sColl^U_\ca$ does have invertible associativity map as well as right unit map. But the left unit map $\lambda$ is not invertible unless \ca is empty (in which case $\sColl^U_\ca$ is just the terminal category with its unique monoidal structure. 

If, as in \cite{BataninMarkl-Operadic}, we adjoint a terminal object to \ca and make this trivial, the resulting operadic category $\ca_1$ will have $R_a\colon \ca/a\to\ca^{|a|}$ exactly as for \ca if $a\in \ca$, while $R_1$ is the identity. Thus $\sColl^U_{\ca_1}$ is in fact monoidal. In particular, if $\ca$ is empty, then $\ca_1$ is the terminal category with its unique operadic structure, and $\sColl^U_{\ca_1}$ is \Set, with the cartesian monoidal structure. 

We now turn to some other examples of operadic categories considered in \cite{BataninMarkl-Operadic}.

\begin{example}
Let \ca be a skeletal abelian category. Choose a representative for
each each quotient. For convenience, suppose that
identity morphisms are chosen quotients. Let $\Epi(\ca)$ be the category whose
objects are those of \ca, and whose morphisms are chosen quotient
maps. Composition is as in \ca, but corrected if necessary to give a
chosen quotient. Define $\Epi(\ca)\to\cs$ to be constant at $1\in
\cs$. Define $R\colon \Epi(\ca)/c\to \Epi(\ca)$ to pick out the kernel (as an object).
Thus $R(r\colon a\to c)=\ker(r)$, and while if $qp=\theta r$ with
$p,q,r$ all chosen quotients, $R(r\colon p\to q)$ is  the restriction
$p^q\colon \ker(r)\to \ker(q)$ of $p$. This defines a genuine  operadic
category \cite[Example~1.22]{BataninMarkl-Operadic}.

A morphism is fibrewise trivial if and only if it is invertible. The functors $R\colon\Epi(\ca)/c\to\Epi(\ca)$ are opfibrations. They are not discrete, but they do reflect isomorphisms by the short five lemma. Thus the 
functors $\theta/R$ are equivalences, and $\sColl^U_{\Epi(\ca)}$ is monoidal.
\end{example}

\begin{example}
 The genuine operadic category $\Bq_I$ of $I$-bouquets was introduced in \cite{BataninMarkl-Operadic} to deal with coloured operads. Here $I$ is a fixed set. An object of $\Bq_I$ is called a bouquet, and consists of an object $m\in\cs$, an $m$-indexed family $c=(c_i)_{i\in m}$ of elements of $I$, and another element $c'\in I$. We write $(m,c,c')$ for such a bouquet. A morphism $(m,c,c')\to (n,d,d')$ can exist only if $c'=d'$, in which case it consists of a morphism $f\colon m\to n$ in \cs. Thus there is an evident forgetful functor $\Bq_I\to\cs$ which defines the cardinality functor. 

For a morphism $f\colon (m,c,c')\to (n,d,c')$ and an element $j\in n$, the corresponding fibre is $(f^{-1}j,c|_{f^{-1}j},d_j)$, where $c|_{f^{-1}j}$ is the $f^{-1}j$-indexed family with $(c_{f^{-1}j})_i=c_i$.

A bouquet $(m,c,c')$ is trivial if and only if $m=1$ and the unique element of $c$ is $c'$ itself. A morphism $f\colon (m,c,c')\to (n,d,c')$ is fibrewise trivial if and only if $f$ is bijective and $d_{fi}=c_i$ for all $i\in m$. Thus the fibrewise trivial morphisms are strictly contained within the isomorphisms.

Since $\Bq_I$ is a genuine operadic category $\sColl^U_{\Bq_I}$ will be monoidal if and only if the associativity map $\alpha$ is invertible. We shall see that this is the case using Proposition~\ref{prop:Hopf}.

Let $(p,e,e')\in\Bq_I$. An object of the slice category $\Bq_I/(p,e,e')$ consists of an object $f\colon m\to p$ of $\cs/p$ equipped with an element $c\in I^m$, while a morphism from $(m,f,c)$ to $(n,g,d)$ is a just a morphism $(m,f)\to(n,g)$ in $\cs/p$. There is a commutative square 
\[ \xymatrix{
\Bq_I/(p,e,e') \ar[r] \ar[d]_R & \cs/p \ar[d]^R \\ \Bq^p_I \ar[r] & \cs^p } \]
where the horizontal arrows are induced by the cardinality functor, and the upper horizontal is a surjective equivalence. 
Given an object $(m,f,c)\in\Bq_I/(p,e,e')$, there is induced a further commutative square
\[ \xymatrix{
(m,f,c)/(\Bq_I/(p,e,e')) \ar[r] \ar[d]_{(m,f,c)/R} & (m,f)/(\cs/p) \ar[d]^{(m,f)/R} \\
Rf/\Bq^p_I \ar[r] & Rf/\cs^p } \]
where both horizontal arrows are surjective equivalences and the right vertical is an equivalence. Thus the left vertical is also an equivalence, and a straightforward calculation shows that the counit can be chosen to be fibrewise trivial. 

Thus $\sColl^U_{\Bq_I}$ is a monoidal category; in fact it is the
usual category of collections for symmetric coloured operads
(symmetric multicategories).
\end{example}

\begin{example}
Another example given in \cite{BataninMarkl-Operadic} is the operadic category $\Omega_2$ of 2-trees. 

An object of $\Omega_2$ is an order-preserving morphism $\partial\colon p_2\to p_1$ between finite ordinals (in other words, a morphism in \cp); we denote such an object by $p$. A morphism $\phi\colon m\to p$ is a commutative square
\[ \xymatrix{
m_2 \ar[d]_{\partial} \ar[r]^{\phi_2} & p_2 \ar[d]^{\partial} \\ m_1 \ar[r]_{\phi_1} & p_1 } \]
where $\phi_1$ is an order-preserving function, and $\phi_2$ is a function which is order-preserving on the fibres of $\partial$; in other words, if $i,j\in m_2$ with $i<j$ and $\partial i=\partial j$, then $\phi_2i\le \phi_2j$. 

The cardinality functor sends $p$ to $|p_2|$ (the underlying set of the ordinal $p_2$, and sends $\phi$ to $\phi_2$. 
Then $R\colon \Omega_2/p\to \Omega^{p_2}_2$ sends 
$\sigma\colon m\to p$ to the family $(m^i_2\to m^i_1)$ where $m^i_2 =
\sigma^{-1}_2i$ and $m^i_1=\sigma^{-1}_1\partial i$. 

A morphism $\phi\colon m\to p$ is fibrewise trivial if and only if
$\phi_2\colon m_2\to p_2$ is bijective and the square defining
$\phi$ is a pullback in \Set. 

We shall see that the sufficient condition of Proposition~\ref{prop:Hopf2} is not satisfied. 
To do this, we identify $n$ with $\{1,\ldots,n\}$, and denote a morphism $n\to m$ by its values. Thus ``$1~1~4$'' denotes the function $f\colon 3\to 4$ with $f(1)=f(2)=1$ and $f(3)=4$.

Let $\theta$ be the morphism
$$\xymatrix{
3 \ar[r]^{1~2~2} \ar[d]_{1~1~2} & 2 \ar[d]^{!} \\ 2 \ar[r]_{!} & 1 }$$
and let $\omega_1$ and $\omega_2$ be 
$$\xymatrix{
1 \ar[r] \ar[d]_{1} & 1 \ar[d]^1 \\
2 \ar[r]_{1~2} & 2 }\quad
\xymatrix{
2 \ar[r] \ar@{=}[d] & 1 \ar[d] \\ 2 \ar[r] & 1 }$$
Suppose that $\theta$ factorized as 
$$\xymatrix{
3 \ar[r]^{\phi} \ar[d]_{1~1~2} & m \ar[r]^{\psi} \ar[d]^{\pi} & 2
\ar[d] \\
2 \ar[r]_{i~j} & n \ar[r] & 1 }$$
then clearly $i=\pi\phi 1=\pi\phi 2$ and $j=\pi\phi 3$.
If also the $\omega$s factorized as 
$$\xymatrix{
\{1\} \ar[r]^{\phi} \ar[d]_{1} & \psi^{-1}1 \ar[d]^\pi
\ar[r]^{\sigma_1} & 1 \ar[d]^{1} \\
2 \ar[r]_{i~j} & n \ar[r]_{\sigma'_1} & 2} \quad
\xymatrix{
\{2,3\} \ar[r]^{\phi} \ar[d]_{1~2} & \psi^{-1}2 \ar[r]^{\sigma_2}
\ar[d]_{\pi} & 1 \ar[d] \\
2 \ar[r]_{i~j} & n \ar[r] & 1 }$$
with the $\sigma$s fibrewise trivial; then in particular $\sigma_2$ is
bijective so that $\phi 2=\phi 3$ and $i=j$. But now
$$1=\sigma'_1 i = \sigma'_1 j = 2$$
gives a contradiction. Thus the sufficient condition does not hold in this example. 
It would be interesting to know whether or not $\sColl^U_{\Omega_2}$ is monoidal.
\end{example}

\begin{example}
A similar calculation shows that the operadic category $\Ord_2$ of 2-ordinals, defined in \cite[Example~1.24]{BataninMarkl-Operadic}, does not satisfy the sufficient condition either.  
\end{example}

The next, final, example does not come from \cite{BataninMarkl-Operadic}; rather, it makes a connection with a result of Andrianopoulos \cite{Andrianopoulos-thesis}.

\begin{example}
Let \cb be an arbitrary category with object-set $B$, and define a cardinality functor $\cb\to\cs$   by $|b|=1$ for all $b$. To give an operadic category structure to \cb is
to gie a functor $R_b\colon \cb/b\to \cb$ for each $b\in B$, satisfying the various conditions. Now a category  \cb and these $R_b$ is precisely the dual of the structure considered in \cite[Section~3.3]{Andrianopoulos-thesis}, characterizing what is needed to make $B$ into a skew monoidale (internal skew monoidal structure) in the monoidal bicategory \Span. 

\end{example}

\bibliographystyle{plain}

\end{document}